\documentclass[onefignum,onetabnum]{siamart220329}

\usepackage{amsfonts}
\usepackage{graphicx}
\usepackage{epstopdf}
\usepackage{algorithmic}
\ifpdf
\DeclareGraphicsExtensions{.eps,.pdf,.png,.jpg}
\else
\DeclareGraphicsExtensions{.eps}
\fi

\newsiamremark{remark}{Remark}
\newsiamremark{hypothesis}{Hypothesis}
\crefname{hypothesis}{Hypothesis}{Hypotheses}
\newsiamthm{claim}{Claim}

\headers{Numerical surgery for mean curvature flow of surfaces}{B.~Kov\'{a}cs}

\title{Numerical surgery \\ for mean curvature flow of surfaces\thanks{Submitted to the editors October 31, 2022, revised \today.}
	}

\author{Bal\'{a}zs Kov\'{a}cs\thanks{
			Institute of Mathematics, University of Paderborn, Warburgerstr.~100, 33098 Paderbron, Germany, (\email{balazs.kovacs@math.upb.de}).
		} 
	   }

\usepackage{xcolor}
\usepackage{graphicx}
\usepackage{amsmath,amssymb,bm}
\usepackage{graphicx,graphics,color}
\usepackage{booktabs,subfig} % for tables
\usepackage{tabularx}
\usepackage{listliketab}

\usepackage{color}

\newcommand\bff{{\mathbf f}}

\newcommand\bfn{{\mathbf n}}

\newcommand\bfu{{\mathbf u}}
\newcommand\bfv{{\mathbf v}}
\newcommand\bfw{{\mathbf w}}
\newcommand\bfx{{\mathbf x}}

\newcommand\bfA{{\mathbf A}}

\newcommand\bfH{{\mathbf H}}

\newcommand\bfM{{\mathbf M}}

\newcommand{\Ga}{\varGamma}

\newcommand{\nbg}{\nabla_{\Ga}}

\newcommand{\mat}{\partial^{\bullet}}

\newcommand{\eps}{\varepsilon}

\newcommand{\pa}{\partial}
\newcommand{\R}{\mathbb{R}}

\def \t {(t)}

\def \to {\rightarrow}
\newcommand{\vphi}{\varphi}

\newcommand{\wtu}{\widetilde{\bfu}}

\newcommand{\wtx}{\widetilde{\bfx}}
\newcommand{\bbk}{\color{black}}
\newcommand{\ebk}{\color{black}}

\renewcommand{\S}{\mathbb{S}}
\newcommand{\nS}{\textnormal{S}}

\ifpdf
\hypersetup{
  pdftitle={Numerical surgery for mean curvature flow of surfaces},
  pdfauthor={B.~Kov\'{a}cs}
}
\fi

\begin{document}

\maketitle

% REQUIRED
\begin{abstract}
	A numerical algorithm for mean curvature flow of closed mean convex surfaces with surgery is proposed. 
	The method uses a finite element based mean curvature flow algorithm based on a coupled partial differential equation system which directly provides an approximation for mean curvature and outward unit normal. The proposed numerical surgery process closely follows the analytical surgery of Huisken \& Sinestrari, and Brendle \& Huisken.
	The numerical surgery approach is described in detail, along with extensions to other geometric flows and methods.
	Numerical experiments report on the performance of the numerical surgery process.
\end{abstract}

% REQUIRED
\begin{keywords}
	geometric flows, mean curvature flow, surgery, singularity, topology change, pinch-off, coupled geometric system, evolving surface finite elements
\end{keywords}

% REQUIRED
\begin{MSCcodes}
	35R01, % PDEs on manifolds
	53C44, % Geometric evolution equations (mean curvature flow, Ricci flow, etc.)
	65M60, % Finite element, Rayleigh-Ritz and Galerkin methods for initial value and initial-boundary value problems involving PDEs
	65M15, % Error bounds for initial value and initial-boundary value problems involving PDEs
	65M12. % Stability and convergence of numerical methods for initial value and initial-boundary value problems involving PDEs
\end{MSCcodes}

\section{Introduction}

The main goal of this paper is to propose an algorithm for \emph{numerical mean curvature flow with surgery} for \emph{closed mean convex} surfaces in $\R^3$. 

This surgery algorithm uses the provably convergent method developed in \cite{MCF}, which is based on discretising a system which includes evolution equations for the mean curvature and the normal vector field, originally derived by Huisken \cite{Huisken1984}. This allows a straightforward adaptation of the surgery approach of Huisken and Sinestrari \cite{HuiskenSinestrari} and Brendle and Huisken \cite{BrendleHuisken2016}. 
To our knowledge, with the exception of the present work, this well-established surgery approach had not been translated into a numerical method.

The numerical surgery approach for mean curvature flow herein is inspired by---and shares---the main ideas of \cite{HuiskenSinestrari,BrendleHuisken2016}:
\begin{itemize}
	\item[-] The flow is stopped slightly before a singularity occurs, i.e.~when the mean curvature $H$ exceeds a threshold.
	\item[-] The surface regions with high curvature are removed, and replaced with more regular ones such that the curvature is decreased.
	\item[-] The flow is then evolved further until another singularity.
\end{itemize}

The method we proposed in \cite{MCF} provides not just an approximation of the surface, but an approximation mean curvature $H_h^n \approx H$, with optimal-order $H^1$-norm error bounds. The fully discrete algorithm of \cite{MCF} is known to be convergent until a singularity occurs, error estimates through the singularities are not shown. This is the key property which allows to follow the steps above. 

\begin{figure}[htbp]
	\centering
	\includegraphics[width=1\textwidth,clip,trim={80 75 90 65}]
	{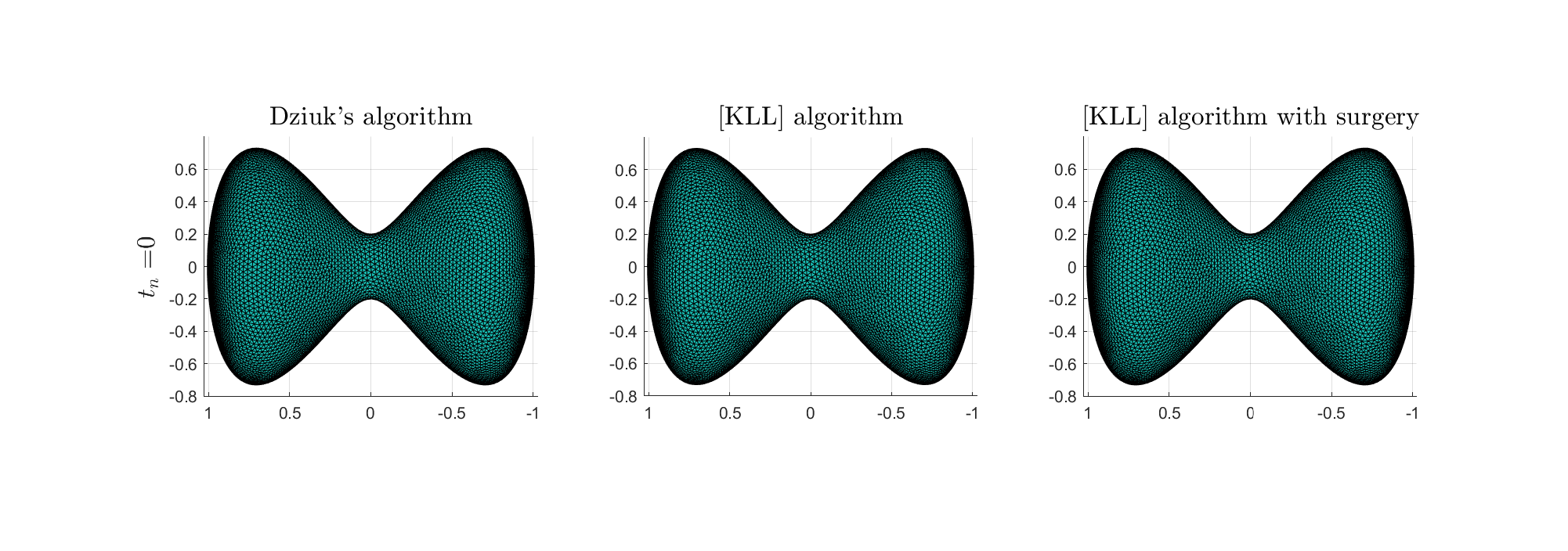}
	\includegraphics[width=1\textwidth,clip,trim={80 75 90 80}]
	{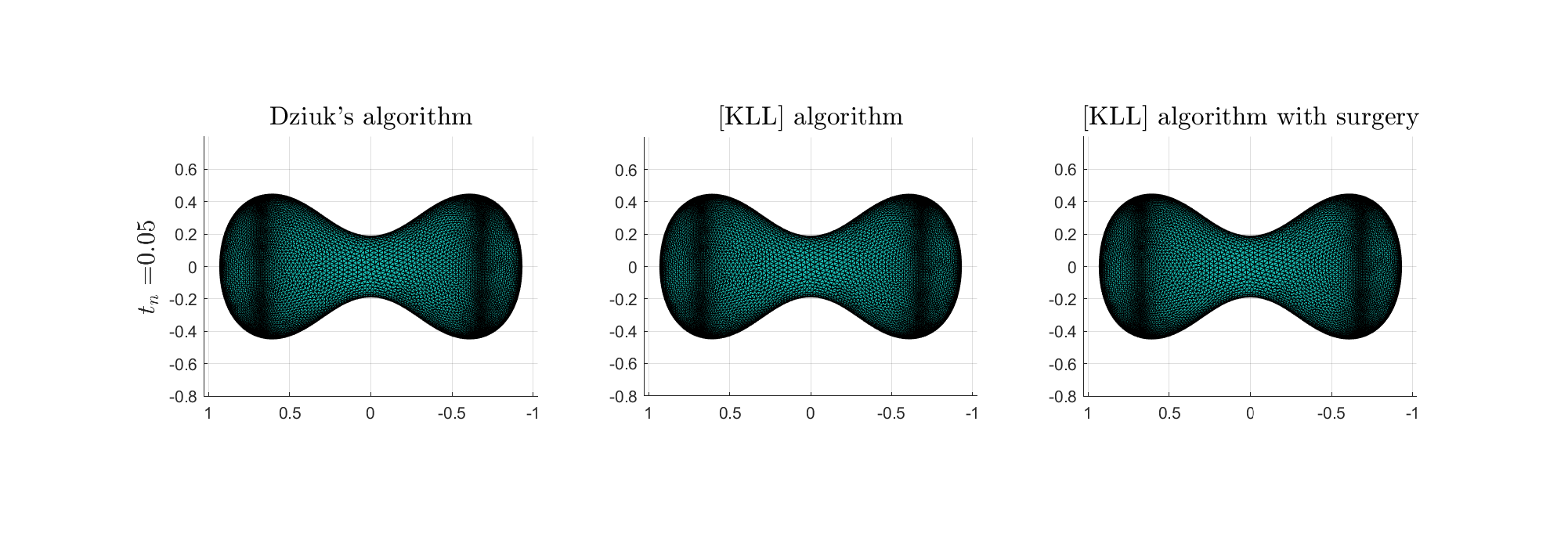}
	% singularity
	\includegraphics[width=1\textwidth,clip,trim={80 75 90 80}]
	{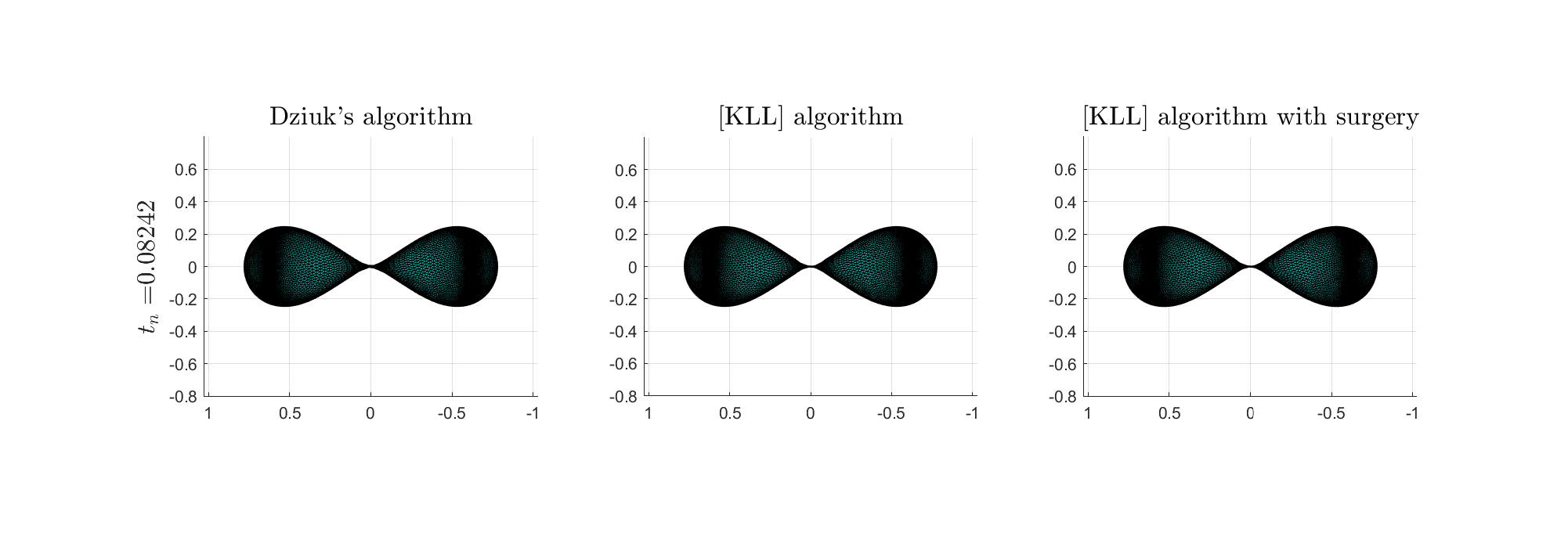}
	\includegraphics[width=1\textwidth,clip,trim={80 75 90 80}]
	{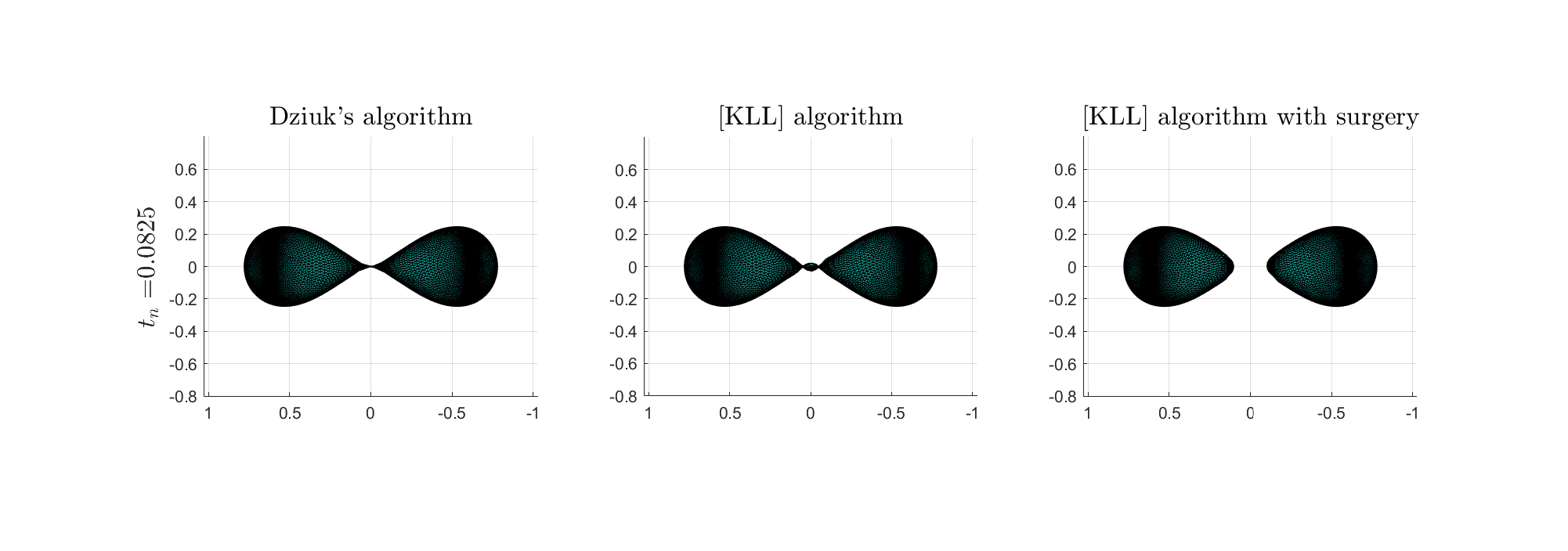}
	\includegraphics[width=1\textwidth,clip,trim={80 75 90 80}]
	{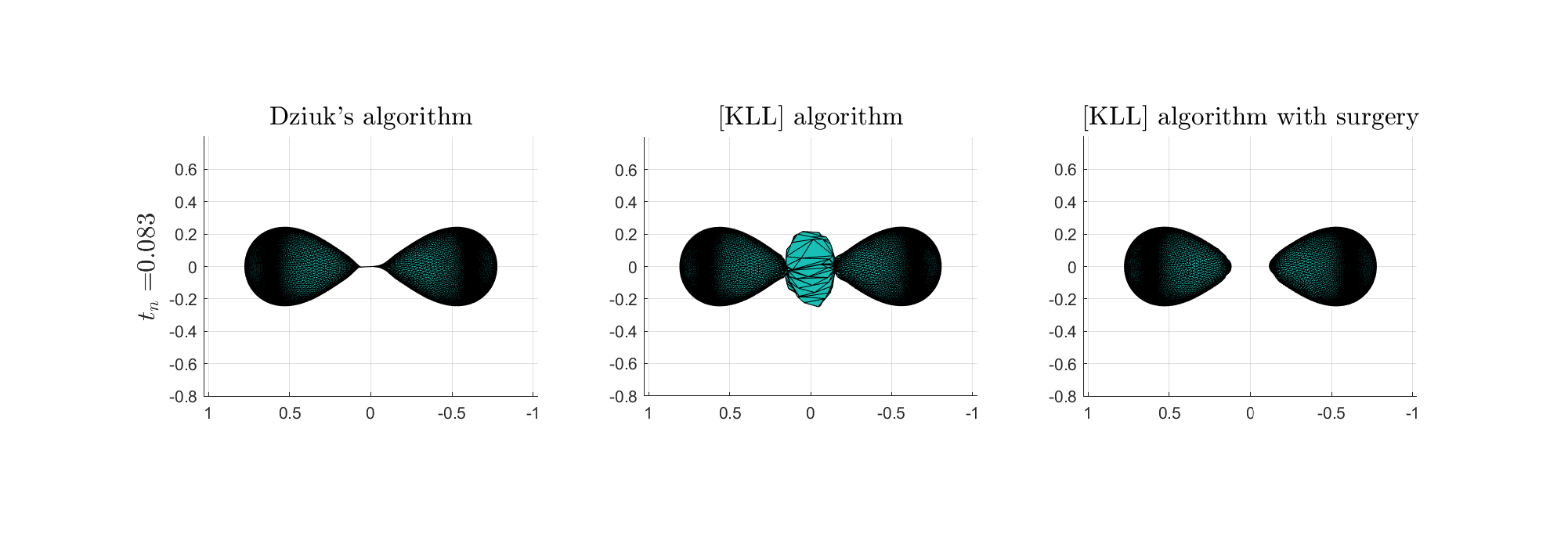}
	\includegraphics[width=1\textwidth,clip,trim={80 75 90 80}]
	{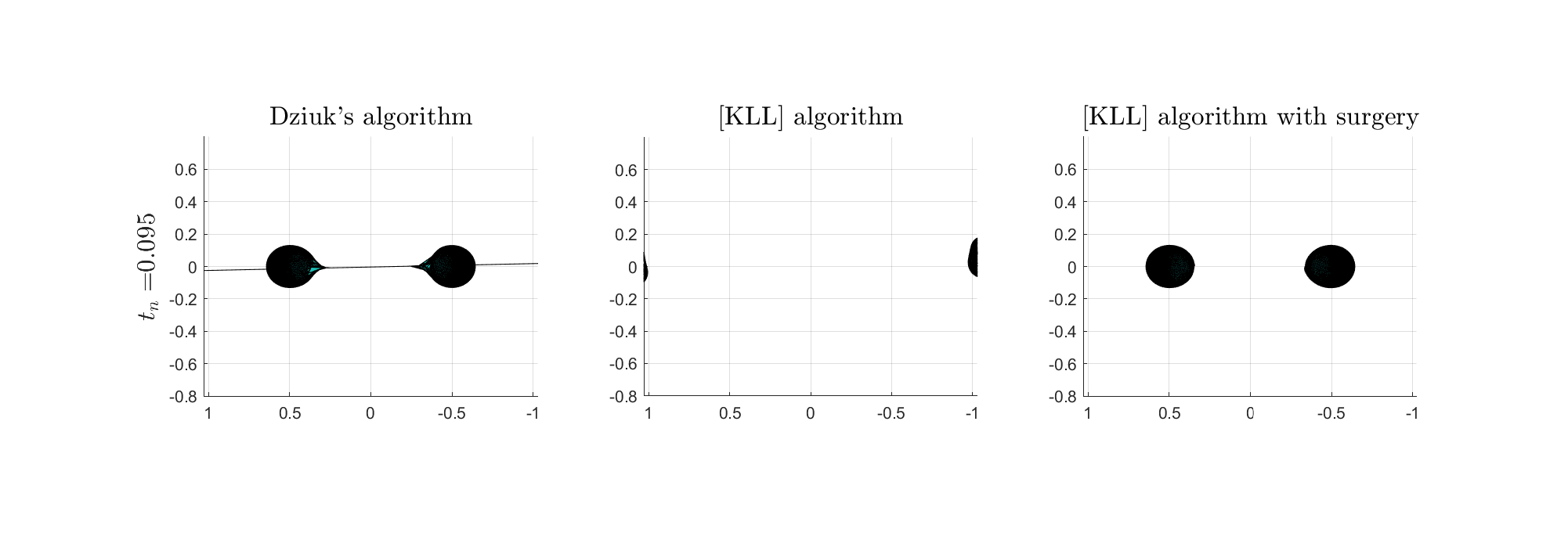}
	\caption{Comparing the numerical solutions of Dziuk's and the [KLL] algorithm  \emph{without surgery} (\cite{Dziuk90} and \cite{MCF}, respectively, on the left and middle columns) with the numerical mean curvature flow  \emph{with surgery} proposed in this paper (right column). A pinch singularity occurs at $t \approx 0.08242$. \newline The surfaces of Dziuk's and the [KLL] algorithm remain connected for all times (note, respectively, the spikes and the overshooting). The [KLL] algorithm with surgery nicely shrinks the two components to round points.}
	\label{fig:compare3}
\end{figure}

The numerical approximation using the proposed approach can be seen in Figure~\ref{fig:compare3}, compared to the methods from \cite{Dziuk90} and \cite{MCF}. More extensive numerical experiments are presented in Section~\ref{section:numerics}.

\medskip
Naturally, it is possible to combine this numerical surgery with other reliable algorithms for mean curvature flow, as long as discrete mean curvature and surface normal are readily available. In particular the convergent method of \cite{HuLi_2022}, which extends the approach of \cite{MCF} and naturally allows tangential motion of the nodes.

We will also give some possible extensions for the numerical surgery of other geometric flows.

\medskip
\textbf{Overview on surgery for geometric flows.} 
Although our work focuses on mean curvature flow, we give some references on the general theory of geometric flows with surgery. 

The formation of singularities in geometric flows has been the target of intensive research since the works of Hamilton \cite{Hamilton95,Hamilton97} for Ricci flow. 
This line of work culminated in the works of Perelman \cite{Perelman_1,Perelman_2,Perelman_3} (using a different surgery process\footnote{Perelman performs surgery \emph{at} singular times, rather then \emph{slightly before} as Hamilton.}) on the Poincaré and geometrisation conjecture.

For mean curvature flow Huisken and Sinestrari translated Hamilton's surgery in \cite{HuiskenSinestrari} for two-convex surfaces of dimension at least $3$. The same result for at least two-dimensional mean convex surfaces was shown, after some refined analysis was available, in \cite{BrendleHuisken2016}; and then in \cite{BrendleHuisken2018} for mean convex surfaces in three-manifolds. An alternative approach was proposed by Haslhofer and Kleiner \cite{HK_2017}.

For more details we refer to \bbk \cite{CMP_survey,HaslhoferKleiner_2017_survey}, and to the references in all these papers. \ebk 

\medskip
\textbf{Overview on numerical methods for mean curvature flow.} 
We give a brief overview on numerical approaches for mean curvature flow of surfaces.

Methods based on surface finite elements (i.e.~parametric methods) go back to the pioneering work of Dziuk \cite{Dziuk90}; Barrett, Garcke, and N\"urnberg used a different velocity-law in \cite{BGN2008}; Elliott and Fritz \cite{Elliott-Fritz-2017} proposed a scheme using the DeTurck trick. The first provably convergent algorithm was proposed in \cite{MCF}. Later on semi-discrete error estimates for Dziuk's method was proved by Li \cite{Li_2021} (and \cite{BaiLi_2022}) requiring finite elements of degree $k \geq 6$. We refer to these papers and their references for more details on parametric methods for mean curvature flow, and also to the review \cite{DeckelnickDE2005} (up to 2005). 
We will give references on finite element based numerical methods for more intricate geometric flows later on.

In \cite{BMPU_2012,MikulaUrban_2012} the authors have proposed a numerical surgery process for curves evolving under forced curve shortening flow. They construct a rectangular background grid which is used to detect (e.g.) pinch-offs. This approach was adapted for a network of curves in \cite{BG_1}, and later to forced mean curvature flow of surfaces in \cite{BG_4}. We are not aware of any other surgery process proposed in the literature.

We also note that, parametric schemes were occasionally integrated beyond a singularity. These experiments have not attempted to resolve these singularities, but merely showed a certain robustness of the algorithms. For instance the first mean curvature flow algorithm of Dziuk already showed such an experiment \cite[Figure~3]{Dziuk90}, further examples are: \cite[Figure~24 \& 25]{BGN2008}, \cite[Figure~4.4--4.6, 4.9, \& 4.12]{BaoJiangZhao_2020}, \cite[Figure~4.12 \& 4.14]{BaoZhao_2021}, % solid-state dewetting
and Figure~\ref{fig:compare3} herein.

Various non-parametric methods, which resolve singularities by design, for mean curvature flow were studied in the literature:
In \cite{Deckelnick2000} Deckelnick proved error estimates for a finite difference scheme approximating level-set (viscosity) solution of mean curvature flow, see as well \cite{DeckelnickDziuk2003}, \cite{OOTT_2011}, and the references therein.
For methods using extended finite elements, see \cite{Olshanskii_etal_2021} approximating the Allen--Cahn equation using trace FEM.
Diffuse interface methods---based on approximation of the Allen--Cahn equation---are, for instance, \cite{FengProhl2003} showing convergence of $O(\eps)$ towards the solution of mean curvature flow, such that only ``little manipulation is needed in order to handle possible singularities''\footnote{Quoted form \cite[Section~1]{FengProhl2003}.}; and see also the results in \cite{Nochetto_etal_1994,Nochetto_etal_1996}. In general, these methods are susceptible to the ``fattening'' of the interface, see, e.g., \cite{FierroPaolini1996}.
Furthermore the thresholding (or MBO) scheme, see \cite{MBO}, also handles singularites very well. 

Non-parametric methods resolve singularities automatically, on the other hand approximation results for these methods are usually limited to error estimates in terms of a tiny (regularisation or phase field) parameter $\eps \ll 1$, or to (weak) convergence without rates, subject to strong CFL-type conditions.

\medskip
\textbf{Outline.}
The paper is structured as follows.
Section~\ref{section:MCF with surgery} gives some details on the surgery of \cite{HuiskenSinestrari} and \cite{BrendleHuisken2016}. 
In Section~\ref{section:numerical MCF with surgery} we first briefly recall our mean curvature flow algorithm from \cite{MCF}, and then describe the main result of this paper: a numerical mean curvature flow with surgery for closed mean convex surfaces. 
In Section~\ref{section:extensions} we give a short but thorough overview on possible extensions for more complicated geometric flows, such as forced mean curvature flow, generalised mean curvature flow, cell-division models, etc.; and we also give details how our numerical surgery process could be used together with other finite element based numerical methods.
In Section~\ref{section:numerics} we present extensive numerical experiments reporting on mean curvature flow through singularity, on the time evolution of maximal mean curvature, and on the proposed numerical surgery process.

\section{Mean curvature flow of surfaces with surgery}
\label{section:MCF with surgery}

Mean curvature flow of a closed sufficiently regular two-dimensional surface $\Ga[X(\cdot,t)] \subset \R^{3}$, with an initial surface $\Ga[X(\cdot,0)] = \Ga^0$, reads:
\begin{equation}
\label{eq:MCF intro}
	\begin{aligned}
		\pa_t X = &\ v \circ X, \\
		v = &\ - H \nu ,
	\end{aligned}
\end{equation}
where $X \colon \Ga^0 \times [0,T] \to \R^3$ is a parametrisation of the surface, $H$ is the mean curvature (i.e.~the sum of principal curvatures $H = \text{tr}(\nbg \nu) = \kappa_1 + \kappa_2$), and $\nu$ is the unit outward normal field. That is, the mean curvature of a sphere is positive.

\subsection{Properties and singularity models for mean curvature flow}
\label{section:MCF properties}

We first recap singularity models and some properties for mean curvature flow, with a particular focus on closed mean convex surfaces.

The parabolic maximum principle implies many important properties for mean curvature flow: preservation of embeddedness, convexity, strict, mean, and of $k$-convexity, and the avoidance property (i.e.~surface do not touch if they started apart, which excludes self-intersections, but allows singularities). See Huisken \cite{Huisken1984}, or, for example, \cite[Section~1.2]{CMP_survey}, \cite[Chapter~3]{Ecker2012}, and \cite[Chapter~2]{Mantegazza}.

The singularities of mean curvature flow of \emph{closed mean convex} hypersurfaces \bbk fall into two classes, type I and type II, separated based on the natural growth rate of the second fundamental form ($A = \nbg \nu$): type I singularity (occuring at time $T$) satisfies the bound
\begin{equation}
\label{eq:natural growth rate}
	\max_{\Ga\t} |A| \leq C \, (T - t)^{-1/2} \qquad \text{for some $C < \infty$,}
\end{equation}
whereas a type II singularity does not follow such a scaling; see~\cite[equation~(1.2)]{AngenentVelazquez_1997}.
The two types \ebk are well characterised:
Huisken and Sinestrari \cite{HuiskenSinestrari_1999a,HuiskenSinestrari_1999b} and White \cite{White_2000,White_2003} have shown that \bbk type I \ebk singularities are either cylindrical or spherical (by \cite{Huisken1990} we even know that the sphere and the cylinder are the only mean convex self-similar shrinkers).
\bbk For type II singularities the large curvature region is convex, differomorphic to a disk, and merges into a cylindrical region of smaller curvature. A detailed description of type II singularities (called the bowl solition) is given in \cite{BrendleChoi_2019,BrendleChoi_2021}, and \cite{Haslhofer_2015}, while a type II singularity example, the degenerate neckpich, was first constructed in \cite{AltschulerAngenentGiga_1995,AngenentVelazquez_1997}. 
We further note that the curvature threshold $H_1$ in \cite{HuiskenSinestrari,BrendleHuisken2016} is used to deal with type II singularities. 
Results on generic solutions of mean curvature flow \cite{MCF_generic_I,MCF_generic_II} support the conjecture that type II singularities are non-generic. \ebk 
For more details (up to 2017) we refer to the surveys \cite{CMP_survey,HaslhoferKleiner_2017_survey}.

%For the surgery process of \emph{two-convex surfaces} ($\kappa_1 + \kappa_2 \geq 0$, \bbk assuming an ordering $\kappa_1 \leq \dotsb \leq \kappa_m$ of the principal curvatures) \ebk $\Ga \subset \R^{m+1}$ of \emph{dimension at least three} ($m \geq 3$) the \emph{cylindrical estimate} \cite[Theorem~1.5]{HuiskenSinestrari} played a crucial role. Roughly, it states that for such surfaces regions of high curvature are either uniformly convex (i.e.~close to a sphere) or close to a cylinder. 
%
%For the surgery of closed embedded \emph{mean convex} ($H > 0$) surfaces $\Ga \subset \R^{3}$ of \emph{dimension two} the cylindrical estimate and its role is replaced by a \emph{sharp inscribed radius estimate} \cite[Theorem~1.1]{Brendle}, which states that, for any $\delta > 0$, the inscribed radius is at least $1 / (1+\delta)H$ at all points where the mean curvature $H$ is greater than $C(\delta)$.

\subsection{Mean curvature flow with surgery}

From \cite{HuiskenSinestrari} and \cite{BrendleHuisken2016} we now recall the surgery construction for mean curvature flow. The former paper by Huisken and Sinestrari established mean curvature flow with surgery for two-convex hypersurfaces in $\R^{m+1}$ of dimension $m \geq 3$, subsequently in the latter paper Brendle and Huisken extended the construction to two-dimensional mean convex surfaces. Their result can be stated as follows.

\begin{theorem}[{\cite[Theorem~1.1]{BrendleHuisken2016}}]
	\label{theorem:BH_surgery} 
	Let $\Ga^0$ be a closed, embedded surface in $\R^{3}$ with positive mean
	curvature. Then there exists a mean curvature flow with surgeries starting from
	$\Ga^0$ which terminates after finitely many steps.
\end{theorem}

Algorithm~\ref{alg:MCF with surgery} briefly describes the main steps of the construction of \cite{HuiskenSinestrari} and \cite{BrendleHuisken2016} for mean curvature flow \emph{with surgery}.

\begin{algorithm}[H]
	\caption{Mean curvature flow with surgery \cite{HuiskenSinestrari} and \cite{BrendleHuisken2016}.}
	\label{alg:MCF with surgery}
	\begin{algorithmic}
		\STATE{\textbf{Data:}~Let $\Ga^0 \subset \R^3$ be a closed embedded initial surface with positive mean \newline \phantom{\textbf{Data:}}~curvature. 
				\newline \phantom{\textbf{Data:}}~Let the curvature thresholds $H_1 < H_2 < H_3$ be given.} \\[1.5mm]
		\STATE{(a) Let the smooth flow \eqref{eq:MCF intro} evolve $\Ga(t)$ until time $t$ such that \\ \phantom{(a)} $\max\{H(\cdot,t)\} > H_3$.}
		\STATE{(b) Perform surgeries on necks, removing all points	with curvature greater than $H_2$. \\ \phantom{(b)} Right after surgery, maximal curvature drops below $H_2$.}
		\STATE{(c) Repeat the Steps (a) and (b) until the flow goes extinct.}
	\end{algorithmic}
\end{algorithm}

Naturally, among the main contributions of \cite{HuiskenSinestrari} and \cite{BrendleHuisken2016} is how exactly, and with what scaling, should one ``perform surgeries on necks'' in Step (b) of Algorithm~\ref{alg:MCF with surgery}. A subtle issue is to show the right $\alpha$-noncollapesedness along the flow, see \cite{BrendleHuisken2016}. For (immersed) surfaces of dimension $m \geq 3$ 2-convexity is required \cite{HuiskenSinestrari,HK_2017}. 
The resulting surfaces have mean curvature that are forced below $H_2$, therefore there holds $\max\{H\} \leq H_2$ except on the regions diffeomorphic to $\S^n$, or to $\S^{n-1} \times \S^1$, which are then discarded. 
Hence, throughout mean curvature flow with surgery the mean curvature $H$ is uniformly bounded by $H_3$.  
Due to its technical nature, for more details (wherein $H_1$ also plays an important role related to type II singularities, see, e.g., \cite[Theorem~2.14 \& 2.15]{BrendleHuisken2016}) we refer to the original works, as well as the reviews \cite{Sinestrari_survey,CMP_survey}.

In \cite{BrendleHuisken2016} the three curvature thresholds are then sent to infinity, showing that the flow with surgery converges to the level set solution. Due to the discrete nature of the numerical computations, we will \emph{not} perform such a limit process (since with a given mesh size one cannot reasonably approximate a surface of arbitrarily large curvature).

We further note that \cite{HuiskenSinestrari} and \cite{BrendleHuisken2016} immediately remove surface components with a trivial topology (e.g.~topological spheres), for the sake of applications, throughout the numerical flow we will keep all components until vanishing.

\section{A numerical algorithm for mean curvature flow with surgery}
\label{section:numerical MCF with surgery}

Before turning to the algorithm for mean curvature flow \emph{with surgery}, we briefly recall the convergent mean curvature flow algorithm from \cite{MCF}, and the corresponding optimal-order error estimates.

\subsection{A convergent algorithm for mean curvature flow until singular times}
\label{section:MCF alg}

The algorithm in \cite{MCF} is based on a system, originally derived by Huisken \cite{Huisken1984}, coupling \eqref{eq:MCF intro} to evolution equations for $H$ and $\nu$ along the flow, and it reads (with $A = \nbg \nu$):
\begin{equation}
\label{eq:MCF system}
\begin{aligned}
v = &\ - H \nu , \\[1mm]
\mat \nu = &\ \varDelta_\Ga \nu + |A|^2 \nu , \\
\mat H = &\ \varDelta_\Ga H + |A|^2 H , \\[1mm]
\pa_t X = &\ v \circ X .
\end{aligned}
\end{equation}
The algorithm discretises the weak form of the above system using evolving surface finite elements (ESFEM) of polynomial degree $k\geq 2$ and linearly implicit backward difference formulae (BDF methods) of order $2$ to $5$ in space and time, respectively. 
The crucial feature of the approach is that all geometric information are computed from the discretised evolution equations, instead of the discrete surface---as it is, for instance, done for the methods in \cite{Dziuk90,BGN2008}, etc.

Without going into deep details of surface finite elements or BDF methods, we consciously recall the algorithm from \cite[Section~2--5]{MCF} (wherein all details can be found), and here only introduce the essential notation.

The surface $\Ga[X]$ and the corresponding variables are approximated in space and time as follows: for a time step size $\tau > 0$ and mesh width $h > 0$, the $N$ basis functions $(\phi_j[\bfx])$, at a time $t_n = n \tau$ we have the approximations
\begin{alignat*}{7}
X(\cdot,t_n) \approx &\ X_h^n \ &=&\ \sum_{j=1}^{N} \bfx^n_j \, \phi_j[\bfx^0] , \qquad\qquad
& v(\cdot,t_n) \approx &\ v_h^n \ &=&\ \sum_{j=1}^{N} \bfv^n_j \, \phi_j[\bfx^n] , \\
\nu(\cdot,t_n) \approx &\ \nu_h^n \ &=&\ \sum_{j=1}^{N} \bfn^n_j \, \phi_j[\bfx^n] , \qquad\qquad
& H(\cdot,t_n) \approx &\ H_h^n \ &=&\ \sum_{j=1}^{N} \bfH^n_j \, \phi_j[\bfx^n] .
\end{alignat*}
That is, the nodal vectors $\bfx^n \in \R^{3N}$, $\bfv^n \in \R^{3N}$, and $\bfn^n \in \R^{3N}$, $\bfH^n \in \R^{N}$---together abbreviated to $\bfu^n = (\bfn^n,\bfH^n)^T$---collect the nodal values of the respective fully discrete approximations at time $t_n$.

These nodal values are determined by the matrix--vector formulation of the fully discrete mean curvature flow algorithm, see~\cite[equation~(5.1)]{MCF}, which reads:
\begin{equation}
\label{eq:BDF}
\begin{aligned}
\bfv^n = &\ - \bfH^n \bullet \bfn^n , \\
\bfM(\wtx^n) \dot \bfu^n + \bfA(\wtx^n) \bfu^n = &\ \bff(\wtx^n,\wtu^n) , \\
\dot \bfx^n = &\ \bfv^n .
\end{aligned}
\end{equation}
The above fully discrete problem uses the stiffness and mass matrices, $\bfA$ and $\bfM$, as well as the nonlinearity $\bff$, the appropriate discretisation of the right-hand side of \eqref{eq:MCF system}. The $\bullet$ denotes the componentwise product of two vectors, see \cite[Section~3.3 \& 5.2]{Willmore}.
The discrete BDF time derivative and the extrapolation are given, for any $(\bfw^n)$, by
\begin{equation*}
\dot \bfw^n = \frac{1}{\tau} \sum_{j=0}^q \delta_j \bfw^{n-j}, \qquad 
\widetilde{\bfw}^n = \sum_{j=0}^{q-1} \gamma_j \bfw^{n-1-j} , \qquad n \geq q ,
\end{equation*}
with coefficients given by $\delta(\zeta)=\sum_{j=0}^q \delta_j \zeta^j=\sum_{\ell=1}^q \frac{1}{\ell}(1-\zeta)^\ell$ and $\gamma(\zeta) = \sum_{j=0}^{q-1} \gamma_j \zeta^j = (1 - (1-\zeta)^q)/\zeta$. %, see\cite[Chapter~V]{HairerWannerII}, \cite{AkrivisLubich_quasilinBDF}.
The method is endowed with initial values, which are computed using a sufficiently accurate method.
For more details we direct the reader, in particular, to Section~3 and 5 in \cite{MCF}.

In \cite{MCF} we proved the following fully discrete error estimate for mean curvature flow, compactly recalled below.
\begin{theorem}[{\cite[Theorem~6.1]{MCF}}]
	\label{theorem:MainTHM-full} 
	The ESFEM--BDF full discretization \eqref{eq:BDF} of the coupled mean curvature flow problem \eqref{eq:MCF system}, using evolving surface finite elements of polynomial degree~$k\geq 2$ and linearly implicit BDF time discretization of order $q$ with $2\le q\le 5$.
	Suppose that the mean curvature flow admits a sufficiently smooth exact solution, that is the flow map $X(\cdot,t):\Ga^0\to \Ga(t)\subset\R^3$ is non-degenerate over $t\in[0,T]$, i.e.~excluding singularities within $[0,T]$. 
	
	Then, there exist for sufficiently small mesh sizes $h > 0$  and time step sizes $\tau >0$ satisfying the step size restriction $\tau \le C_0 h $, 
	the following error bounds hold for the lifts of the discrete position and mean curvature
	\begin{align*}
	\|(\mathrm{id}_{\Ga[X_h^n]})^L - \mathrm{id}_{\Ga(t_n)}\|_{H^1(\Ga(t_n))^3} &\leq C(h^k+\tau^q), \\
	\|(H_h^n)^L - H(\cdot,t_n)\|_{H^1(\Ga(t_n))} &\leq C(h^k+\tau^q) . 
	\end{align*}
	Analogous optimal-order $H^1$ norm error estimates hold for velocity $v$, and normal vector $\nu$.
\end{theorem}

The regularity assumptions required for Theorem~\ref{theorem:MainTHM-full} exclude the formation of singularities, therefore convergence through surgery is not shown here. 

\subsection{Fully discrete mean curvature flow with surgery}

The main goal of this work is to extend the approach of Algorithm~\ref{alg:MCF with surgery} (under the same assumptions, in particular including mean convexity) to the parametric numerical approach of Section~\ref{section:MCF alg} above.

A crucial property---especially with regard to the goal of this paper---of the method \eqref{eq:BDF} from \cite{MCF} is that we have direct access to an approximation of the mean curvature $H_h^n$. In addition, the optimal-order error estimates provided by Theorem~\ref{theorem:MainTHM-full} make it possible to directly translate the approach of \cite{HuiskenSinestrari,BrendleHuisken2016} to the fully discrete setting of \eqref{eq:BDF} (note that $H_1$ does not play a role here). This almost verbatim translation is found in Algorithm~\ref{alg:numerical MCF with surgery}.

\begin{algorithm}[H]
	\caption{Fully discrete mean curvature flow with surgery.}
	\label{alg:numerical MCF with surgery}
	\begin{algorithmic}
		
		\STATE{\textbf{Data:}~Let $\Ga^0$ be a closed embedded initial surface with positive mean curvature. \newline 
			\phantom{\textbf{Data:}}~Let $\Ga_h^0$ be the interpolated discrete initial surface, and given interpolated \newline 
			\phantom{\textbf{Data:}}~initial data $H_h^0$ and $\nu_h^0$. \newline 
			\phantom{\textbf{Data:}}~Let the curvature thresholds $H_3 > H_2$ be given.
		} \\[1.5mm]
		\STATE{(a) Let the fully discrete algorithm \eqref{eq:BDF} evolve $\Ga_h[\bfx^n]$, $H_h^n$, and $\nu_h^n$ until time $t_n$ \\ \phantom{(a)} such that $\max\{H_h^n\} > H_3$.}
		\STATE{(b) Perform surgeries on $\Ga_h[\bfx^n]$ (see Algorithm~\ref{alg:numerical surgery}): remove all nodes for \\ \phantom{(b)} which $H_h^n(\bfx^n_j) > H_2$, close the surfaces, and compute geometric quantities. \\ \phantom{(b)}
			Right after surgery maximal curvature drops to $H_2$.}
		\STATE{(c) Repeat the Steps (a) and (b) until the discrete flow goes extinct.}
	\end{algorithmic}
\end{algorithm}

There are only minimal differences between the continuous and numerical approaches, compare Algorithm~\ref{alg:MCF with surgery} and \ref{alg:numerical MCF with surgery}. 
The main part of the numerical approach, Step (b), is described in detail in Section~\ref{section:surgery}

\begin{remark}
\label{remark:other parametric methods}
	For other parametric methods, e.g.~\cite{Dziuk90,BGN2008}, such a verbatim translation would be far more problematic, since these methods for mean curvature flow directly compute geometric informations form the discrete surface. 
	Therefore a numerical surgery for the numerical flows of these methods would need to detect singularities. A possibility is to detect close points which are far along a geodesic (e.g., see below). Or by computing the missing geometric information directly from $\Ga_h$, similarly to Step (d) in Algorithm~\ref{alg:numerical surgery}. In either cases, contrary to the method \eqref{eq:BDF} from \cite{MCF}, the question of accuracy (i.e.~error estimates), for both the numerical method and the geometric computations, remains an open issue. 
	
	As mentioned in the introduction, for curves evolving under forced curve shortening flow, a background grid based approach was proposed in \cite{BMPU_2012,MikulaUrban_2012}, and further explored and extended for a network of curves in \cite{BG_1}, and \cite{BG_4} extended this approach to surfaces in three-dimensions. Note that this background grid method can also handle coalescence, and triple-junctions.
	
	Since mean curvature flow does not develop self-intersections, see Section~\ref{section:MCF properties}, another possibility to detect singularities would be to define (and evaluate in each time step) a non-local energy which penalises self-intersections, whence high energy would indicate a singularity. For self-avoiding curves \cite{BRR_2018} used a ``repulsive tangent-point potential'' proposed in \cite{GonzalezMaddocks_1999}.
\end{remark}

\subsubsection{An algorithm for numerical surgery}
\label{section:surgery}

We will now describe the numerical surgery process, Step (b) in Algorithm~\ref{alg:numerical MCF with surgery}. 

Regarding terminology: under a triangle or an element we will mean an element of polynomial degree $k$ of the discrete surface $\Ga_h$, i.e.~the image of the unit simplex under a degree $k$ reference mapping, see~\cite{Demlow2009,highorderESFEM,MCF}. 

Assume at some time $t_n$ the condition $\max\{H_h^n\} > H_3$ is true, and hence the numerical method \eqref{eq:BDF} is stopped momentarily, and a numerical surgery (analogous to Algorithm~\eqref{alg:MCF with surgery}.b) shall be performed. 

The numerical surgery has four main steps summarised in Algorithm~\ref{alg:numerical surgery} below.

\begin{algorithm}[H]
	\caption{Numerical surgery algorithm.}
	\label{alg:numerical surgery}
	\begin{algorithmic}
			\STATE{\textbf{Data:}~Pre-surgery surface $\Ga_h^n := \Ga_h[\bfx^n]$ given by its nodes and elements, and \newline \phantom{\textbf{Data:}}~its discrete geometric information $\nu_h^n$ and $H_h^n$. \newline 
				\phantom{\textbf{Data:}}~A set of nodes $\{\bfx^n_j\}$ marked for removal by the rule $H_h^n(\bfx^n_j) > H_2$.
			}
			\STATE{\textbf{Result:}~Post-surgery surface $\Ga_h^{n+} := \Ga_h^\nS$ given by its nodes and elements, and the \newline \phantom{\textbf{Result:}}~corresponding geometric information $\nu_{\Ga_h^\S}$ and $H_{\Ga_h^\S}$.
			} \\[1.5mm]
			\STATE{(a) Remove all elements of $\Ga_h^n$ which have a node satisfying $\max\{H_h^n\} > H_2$. \\ \phantom{(a)} Remove all nodes of the removed elements.}
			\STATE{(b) Determine the boundaries of the remaining surface parts. \\ \phantom{(b)} Based on the connectivity provided by the elements of $\Ga_h^n$ determine the \\ \phantom{(b)} topology and orientation of these boundaries.}
			\STATE{(c) Close all boundaries of the obtained surface parts:}
			\FOR{\textnormal{all boundaries}}
				\STATE{(c.0) Generate a spherical cap which has the same number of edges as \\ \phantom{(c.0)} the boundary.}
				
				\STATE{(c.1) Determine the position, height and radius of the cap, such that its curvature \\ 
					\phantom{(c.1)} is not exceeding $H_2$.}
				
				\STATE{(c.2) Sew together the cap and the part, by adding nodes and elements.}
			\ENDFOR
			\STATE{This yields the post-surgery surface $\Ga_h^\nS$.} \\[1.5mm]
			\STATE{(d) Compute the necessary geometric quantities $\nu_{\Ga_h^\nS}$ and $H_{\Ga_h^\nS}$,  if they cannot be \\ \phantom{(d)} retained from $\nu_h^n$ and $H_h^n$, or obtained directly from the spherical cap.}
	\end{algorithmic}
\end{algorithm}

\bigskip
We would like to give further details on Step (c.0)--(c.2), and (d) of Algorithm~\ref{alg:numerical surgery}.

\textit{Step (c.0)~~Generating spherical caps.} 
The spherical caps (hemi-spheres of unit radius) \bbk are generated \ebk by an iterative process to ensure that the cap and the part have the same number of edges. 
For the surface mesh generation a modified version of DistMesh \cite{distmesh} is used which handles surfaces with boundary. 
A crude iterative process is used to generate an initial cap which has almost as many edges as the current boundary. Then some of the elements with longest boundary edges are bisected until boundary lengths match.
Finally, we apply a few DistMesh steps equidistributing the nodes on the boundary and cap surface.

\textit{Step (c.1)~~Cap size and position.} 
First we determine the size of the cap. Using the boundary nodes of the current part, we compute the radius of the base ($R_{\text{b}}$) of the cap its height, and the radius of the sphere ($R_{\text{S}}$) it lies in, such that there is still a gap between the boundary and the cap base which is used to connect the two surface with elements proportional to the smallest boundary edge. 
Then the cap is positioned (rotated and translated) to its final position such that its points into the same direction as the vector obtained from the orientation of the removed thin cylindrical part.

\textit{Step (c.2)~~Sewing.} 
While sewing together two surfaces we walk through all boundary edges of the part and the cap, starting from a pair of opposite nodes which are closest to each other. Then, each edge pair is attached by two elements forming a quadrilateral between the two edges. From the two diagonals the shorter one is chosen, to preserve mesh regularity as much as possible.
Since our numerical algorithm \eqref{eq:BDF} uses at least quadratic finite elements, new nodes are added which are linear combination of opposing edge points.
%This process is shown in Figure~\ref{fig:sewing}.

This explains the requirement in Step~(c.0), that a spherical cap should be generated that has the same number of edges as the boundary. In fact, the spherical cap is only generated with a given boundary edge count, if it was not stored during a previous computation before.

\begin{figure}[htbp]
	\centering
	\includegraphics[width=0.48\linewidth,clip,trim={290 60 270 53}]{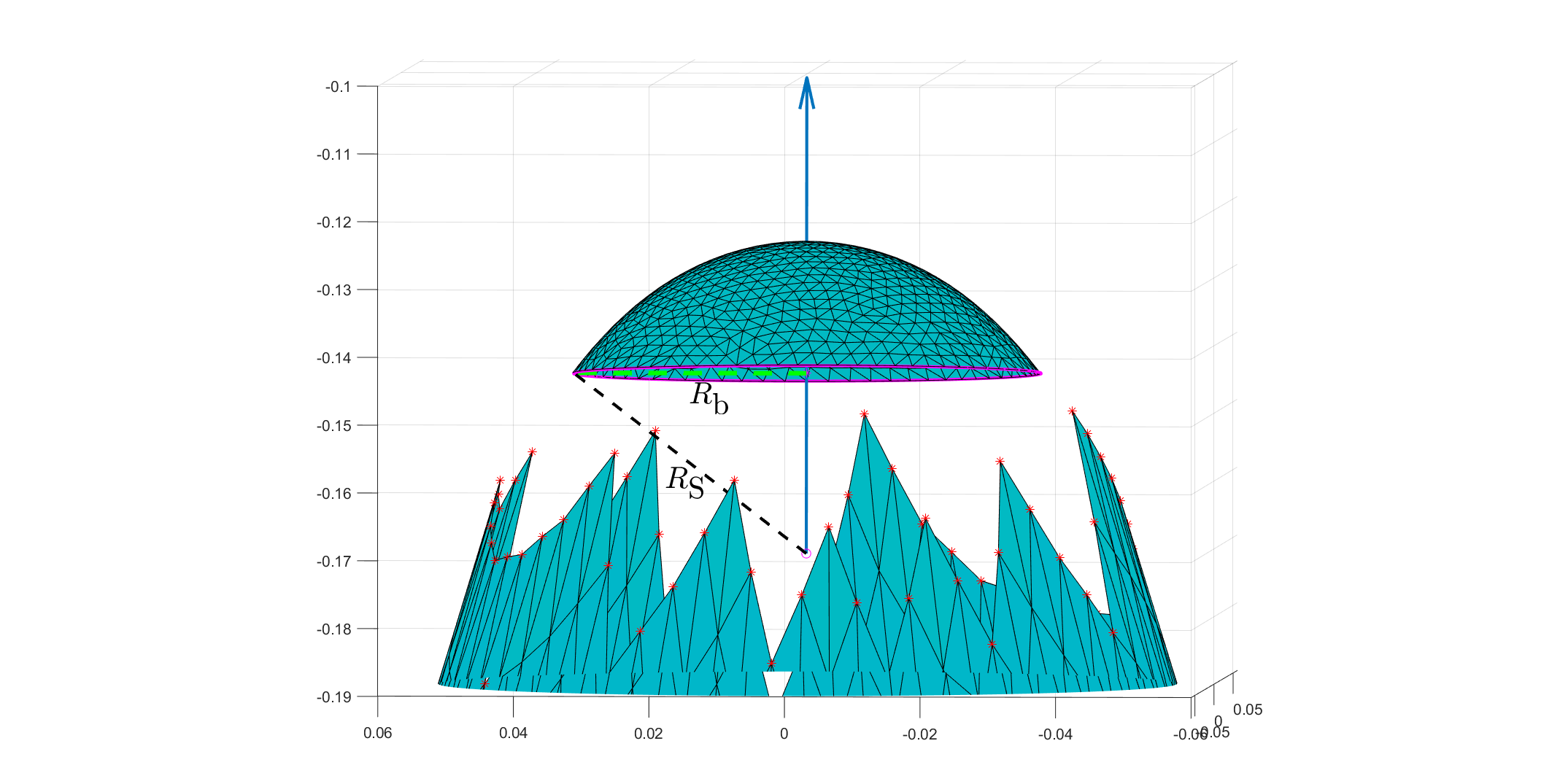}
	\includegraphics[width=0.48\linewidth,clip,trim={250 100 240 100}]{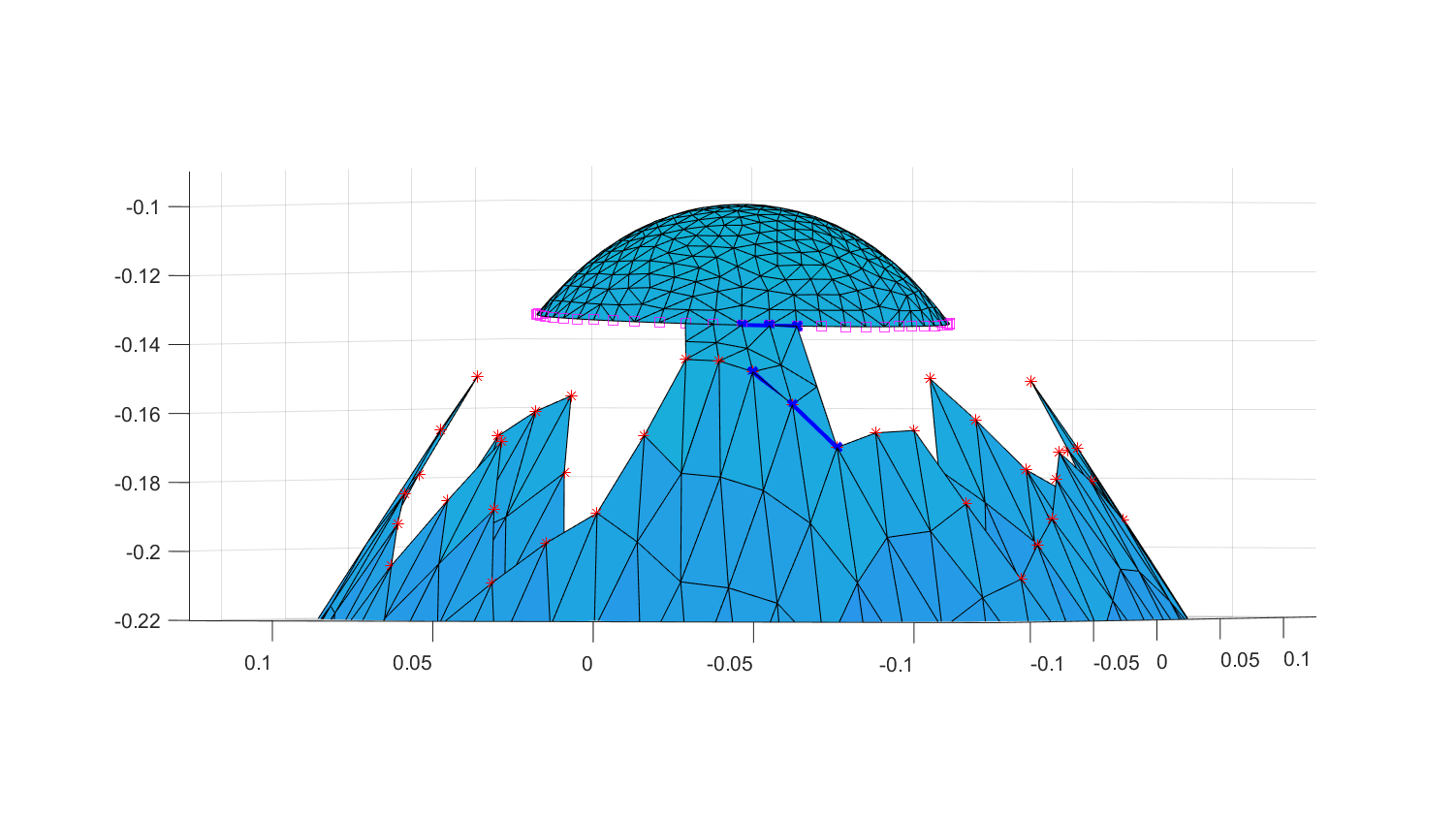}
	\caption{(a) Step (c.1): The scaled and positioned spherical cap, showing the sphere radius $R_\text{S}$ and the base radius $R_\text{b}$, and the axis of the removed cylindrical part. (b) A snapshot of Step (c.2): The current edges (thick lines) on the two boundaries (marked by stars and squares) are connected with two new elements spanned between them.}
	\label{fig:sewing}
\end{figure}

\textit{Step (d): Geometric quantities $\nu_{\Ga_h^\nS}$ and $H_{\Ga_h^\nS}$.} 
The geometric quantities are only computed if they cannot be retained either from the pre-surgery surface $\Ga_h^n$, or obtained from the spherical cap.

The normal vector $\nu_{\Ga_h^\nS}$ is directly computed from the elements of the discrete surface $\Ga_h^\nS$, by taking the area-weighted linear combination of the normal vectors corresponding to a patch, exactly as in~\cite[equation (2.7)]{BGN2008}.

To compute $H_{\Ga_h^\nS}$, motivated by the identity $H = \textnormal{div}_\Ga \, \nu$, we solve the discretisation of the elliptic problem for $\widehat \Ga \subset \Ga$:
\begin{equation*}
\int_{\widehat \Ga} H \vphi = - \int_{\widehat \Ga} \nu_{\Ga} \cdot \nbg \vphi , \qquad (\vphi \in H_0^1(\widehat \Ga) ) ,
\end{equation*}
endowed with inhomogeneous Dirichlet boundary conditions.

\section{Extensions}
\label{section:extensions}

We will give a brief overview on possible extensions to geometric flows beyond mean convex mean curvature flow, and also discuss numerical methods which are compatible with the surgery approach of Algorithm~\ref{alg:numerical MCF with surgery}--\ref{alg:numerical surgery}.

\subsection{Applications}

There are mathematical models describing various real-life phenomena for which geometric flows play a central role. 
Below, if possible, we immediately give references to surgery-compatible algorithms.

These models, instead of mean curvature flow $v = - H \nu$, set the velocity law $v = V \nu$, with normal velocity $V$, to:

\storestyleof{description}
\begin{listliketab}
	\begin{tabularx}{\linewidth}{l X}
		$V = - H + g(u)$ & Forced mean curvature flow where $u$ is the solution of a PDE on the surface, see \cite{MCF_soldriven}. This includes many interesting applications: tumour growth \cite{BarreiraElliottMadzvamuse2011}; surface dissolution \cite{EilksElliott2008}, 
		diffusion induced grain boundary motion \cite{FifeCahnElliott2001}, image segmentation \cite{BG_1,BG_4}, etc. See also the references in these papers. We note that self-intersections cannot be directly detected by Algorithm~\ref{alg:numerical MCF with surgery}, see Remark~\ref{remark:self-intersections}. \\
		%\end{tabularx}
		%\end{listliketab}
		%
		%\begin{listliketab}
		%\begin{tabularx}{\linewidth}{l X}
		$V=-V(H)$ & Generalised mean curvature flow \cite{MCF_generalised}, which includes mean and inverse mean curvature flow \cite{Huisken1984,HuiskenPolden,HuiskenIlmanen}, as well as their powers \cite{Schulze_diss,Schulze_1,Gerhardt_pIMCF,Scheuer_pIMCF}, and also the non-homogeneous mean curvature flow of \cite{AlessandroniSinestrari_nhMCF,Espin_nhMCF}, etc. \\ 
		$V=-F(H,u)$ & Mean curvature flow interacting with diffusion \cite{MCFdiff}, which includes the multiplicative forcing $V = -g(u)H$ of \cite{diss_Buerger,ABG_MCFdiff_1,ABG_MCFdiff_2}, and many other flows. \\ 
		$V=\varDelta_\Ga H + Q$ & Fourth-order flows, i.e.~Willmore flow, or surface diffusion, etc., see \cite{Willmore}, and their forced versions, for example, solid-state dewetting \cite{BaoJiangZhao_2020,BaoZhao_2021}.
	\end{tabularx}
\end{listliketab}

%\begin{remark}[Topology of surgical parts]
%	content...
%\end{remark}
\begin{remark}[Self-intersections]
	\label{remark:self-intersections}
	Mean curvature flow does not develop self-inter\-sections (i.e.~the surface merging with itself without the mean curvature blowing up), cf.~Section~\ref{section:MCF properties}. 
	For generalised geometric flows \emph{self-intersections} may occur, for instance see \cite{MCF_soldriven} and the references therein. Self-intersections can even occur as the true solution of the geometric flow, for example in the model \cite[Section~5.3]{diss_Buerger} and \cite[Section~10.5]{MCFdiff}.
	
	If a self-intersection occurs \emph{which need to be resolved} and a geometric variable indicates this (like $H$ did with singularities for mean curvature flow), the proposed approach can be adapted to cover these cases. If none of the variables provide reliable information on self-intersections, the above repulsive tangent-point potential  \cite{GonzalezMaddocks_1999,BRR_2018}, or background grid approach \cite{BMPU_2012,MikulaUrban_2012,BG_1,BG_4} could provide a remedy. 
\end{remark}
\begin{remark}[Non-mean convex surfaces]
\label{remark:fat torus}
	Generalising the numerical surgery Algorithm~\ref{alg:numerical surgery} to the non-mean convex case ($H \ngtr 0$) is straightforward, on the other hand this (currently) requires the \emph{assumption} that the singularities are still restricted to cylinders or spheres. This would cover other singularity types, such as when a torus merges into a sphere.
	Other types of singularities, if characterised analytically, could also be added to the numerical method.
	
	The stopping and node removing criteria should then be modified accordingly, that is Algorithm~\ref{alg:numerical MCF with surgery} is stopped at time $t_n$ such that $\max\{|H_h^n|\} > H_3$, while the numerical surgery Algorithm~\ref{alg:numerical surgery}~(a) removes all nodes which satisfy $|H_h^n(\mathbf{x}_j^n)| > H_2$. In Step~(d), the computed geometry for the generated caps and newly added nodes reflects whether for the removed parts either $H_h^n > H_2$ or $H_h^n < -H_2$ is true:  $\pm \nu_{\Ga_h^\nS}$ and $\pm H_{\Ga_h^\nS}$ in the two cases, respectively. 
	
	A numerical experiment for a torus merging into a sphere is presented in Section~\ref{section:fat torus}.
\end{remark}

\bigskip
There are interesting models which do not fit into the above problem class ($v = V \nu$).

For example, the system derived and numerically studied by Wittwer and Aland \cite{WittwerAland2022} models self-organized cell division through active surfaces in fluids by requiring the velocity $v$ and the surface stress tensor $S_\Ga$ to satisfy
\begin{equation*}
S_\Ga = (1 - \eta) \nbg \cdot v + Q ,
\end{equation*}
where $\eta$ is a given modelling constant, and $Q$ collects some further contributions, see \cite[equation~(7)--(12)]{WittwerAland2022}.

Their simulations, see \cite[Figure~4]{WittwerAland2022}, stop before a cell division indeed occurs. The local concentration of stress generating surface molecule \cite[equation~(6)]{WittwerAland2022} could serves as an indicator for singularities.

\subsection{Numerical methods}

We give here a brief overview of finite element based numerical methods for geometric flows which could be directly combined with the numerical surgery process herein. 

The convergent algorithms proposed and studied in \cite{MCF_soldriven,Willmore,MCF_generalised,MCFdiff} are all based on the discretisation of a coupled geometric system, similar to \eqref{eq:MCF system}, and hence are directly suitable for our surgery process.

The methods of Barrett, Garcke, and N\"urnberg \cite{BGN2008} are based on the velocity law $$v \cdot \nu = V,$$ which allows for natural tangential motion of the nodes. For mean curvature flow, their system computes an approximation of mean curvature, see \cite[equation~(2.4)]{BGN2008}. Therefore, combination with the proposed numerical surgery is possible with the use of this curvature approximation. However, no convergence results are available for the method of \cite{BGN2008}.

Similarly, to the works of Barrett, Garcke, and N\"urnberg the methods of \cite{HuLi_2022} (for mean curvature and Willmore flow) also use the velocity law $v \cdot \nu = V$, and---extending the ideas of \cite{MCF}---couple this to evolution equations for $H$ and $\nu$. This combination allows for tangential motions \emph{and} a proof for optimal-order error estimates. 
%These properties are also ideal for our numerical surgery process.

\section{Numerical experiments}
\label{section:numerics}

We performed numerical simulations and experiments for mean curvature flow with surgery as proposed in Algorithm~\ref{alg:numerical MCF with surgery} and \ref{alg:numerical surgery}, and also for the case described in Remark~\ref{remark:fat torus}.

The numerical experiments use quadratic evolving surface finite elements. The presented figures all display the linear interpolation of the quadratic discrete surface (i.e.~four triangles in the figures compose one element for the ESFEM approximation).
For the matrix assembly the reference mapping of the quadratic elements was inspired by \cite{BCH2006}. We employed quadratures of sufficiently high order to compute the finite element vectors and matrices so that the resulting quadrature errors are sufficiently small. 
We use the 2-step linearly implicit BDF time discretisation, see \eqref{eq:BDF}. Initial values were obtained by a linearly implicit backward Euler step (which has local order two).

Initial meshes were all generated using DistMesh \cite{distmesh}, without taking advantage of any symmetry of the surfaces. Note that the initial values are not necessarily mean convex, this merely saved time finding suitable initial surfaces undergoing singularities.

With a given mesh size $h$ one cannot reasonably approximate a surface of arbitrarily large curvature, therefore the curvature thresholds $H_j$ are needed to be chosen depending on the mesh width (in our experiments, a reasonable choice was $h \cdot H_3 \approx \mathcal{O}(1)$). 
For the same reason, we do not attempt the limit process $H_j \to \infty$. 

More in depth experiments on the choice of curvature thresholds and the limit process should probably be done using a space--time adaptive algorithm, cf.~\cite{DemlowDziuk_2007,Lantelme}. For a given threshold $H_3$, adaptivity would allow to resolve a singularity much better in both space and time.

% Dumbbell: dof = 10.5k; 0.04899 * 200
% TorusSphere: dof = 11.4k; 0.30815 * 20
% Torus_fatG: dof = 11.4k; (0.7145 , 0.0326) * 200

\subsection{Numerical mean curvature flow with surgery for a dumbbell}

For this experiment the initial dumbbell-shaped surface---which develops a pinch singula\-ri\-ty---is given as $\Ga^0 = \{x \in \R^3 \mid d(x) = 0 \}$ with the signed distance function
\begin{equation}
\label{eq:distance - dumbbell}
d(x) = x_1^2 + x_2^2 + 2 x_3^2 \big( x_3^2-199/200 \big) - 0.04 .
\end{equation}

\begin{figure}[htbp]
	\centering
	\includegraphics[width=1\textwidth,clip,trim={82 58 50 7}]
	{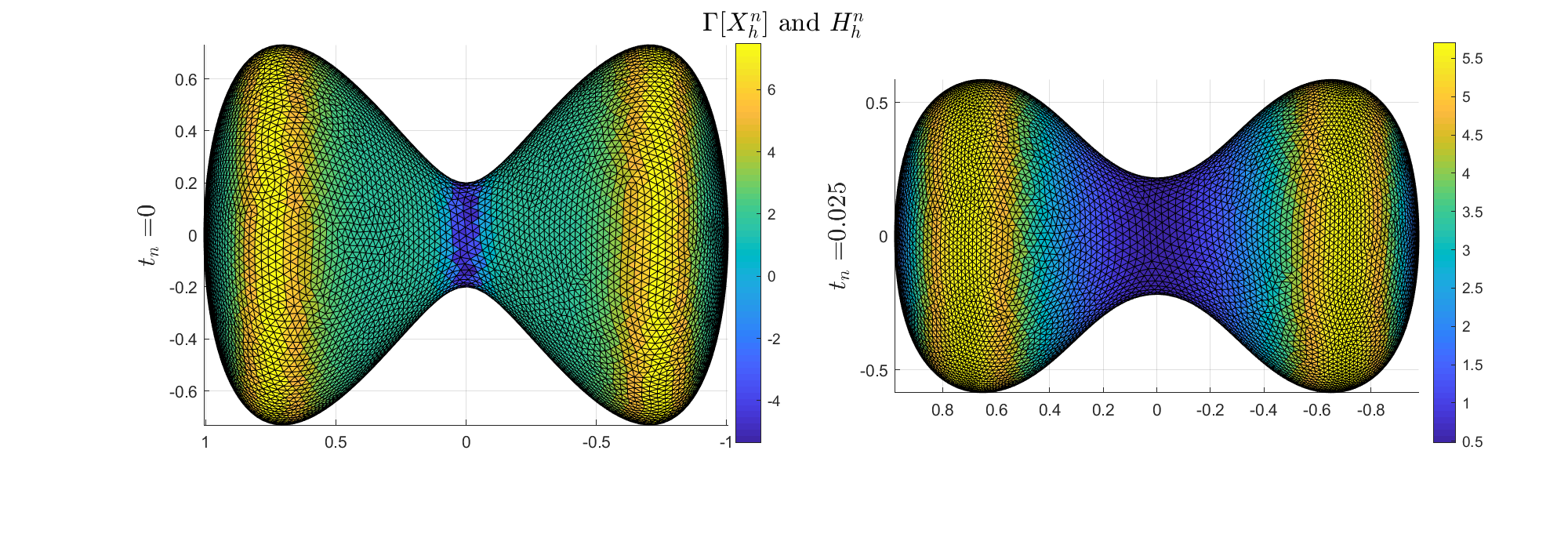}
	\includegraphics[width=1\textwidth,clip,trim={82 80 50 70}]
	{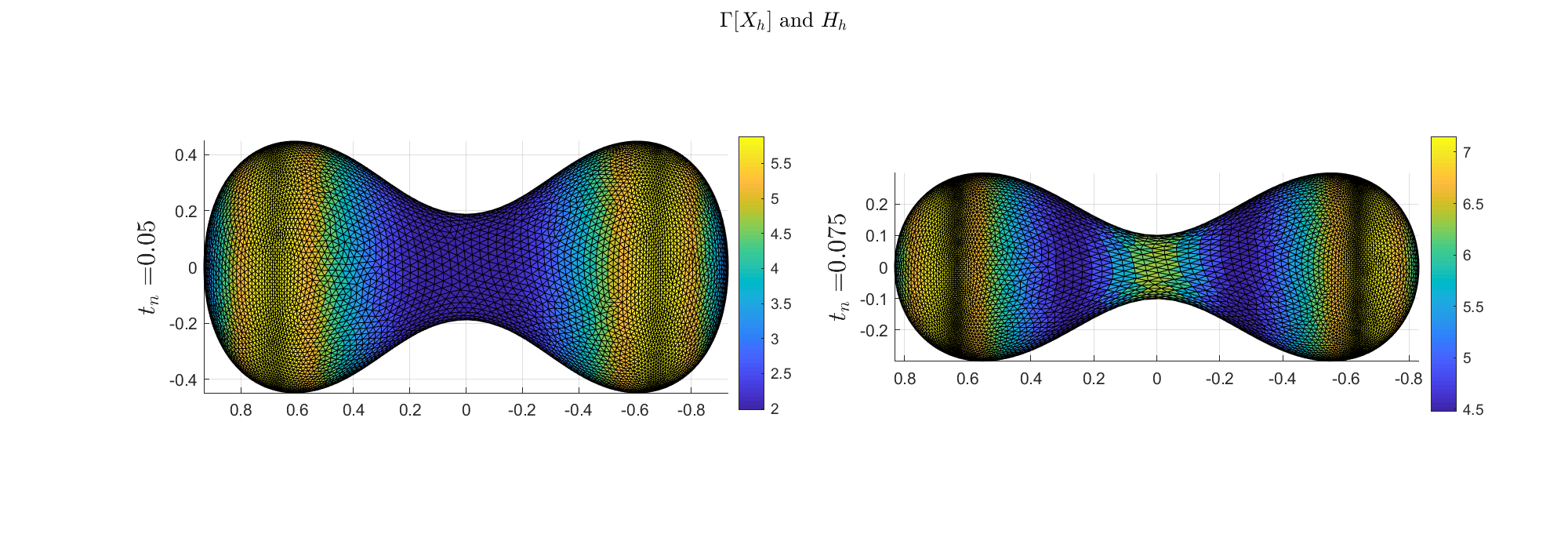}
	\includegraphics[width=1\textwidth,clip,trim={82 108 47 85}]
	{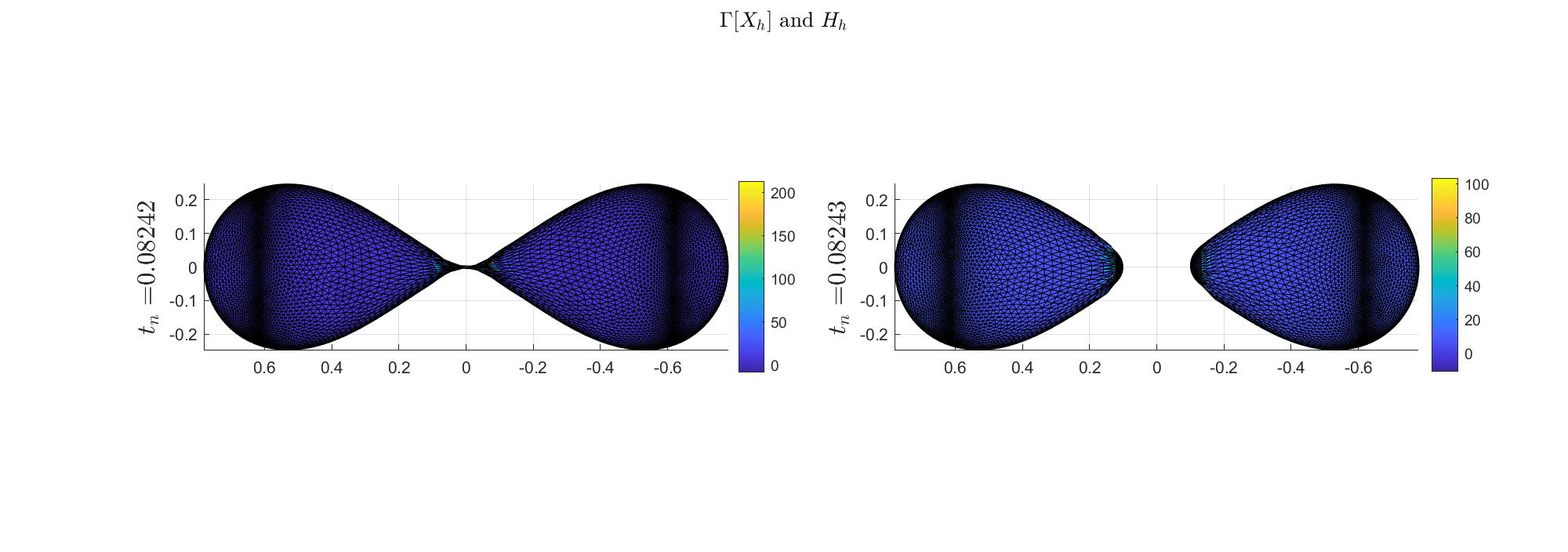}
	\includegraphics[width=1\textwidth,clip,trim={82 115 47 90}]
	{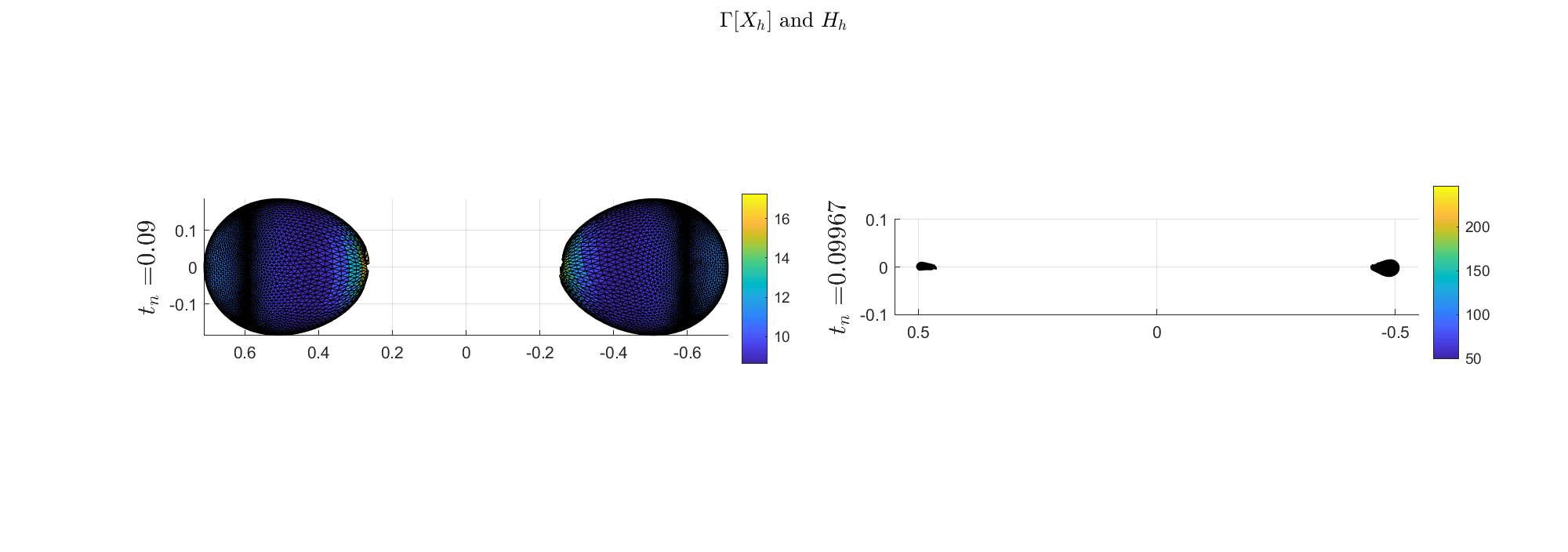}
	\caption{The numerical solution, $\Ga[X_h^n]$ (rotated), and $H_h^n$ at various time steps $t_n$, of mean curvature flow with surgery for a dumbbell. Surgery is plotted in the third row.}
	\label{fig:MCF_with_surgery}
\end{figure}

Figure~\ref{fig:MCF_with_surgery} reports on numerical mean curvature flow with surgery for this dumbbell surface.\footnote{A video of this experiment is available at \href{https://kovacs.app.uni-regensburg.de/research.html}{\texttt{https://kovacs.app.uni-regensburg.de/research.html}}.} The interpolation of the initial surface $\Ga^0$ (with \eqref{eq:distance - dumbbell}, see top left in Figure~\ref{fig:MCF_with_surgery}; $\text{dof} = 10522$ and $h = 0.04899$) is evolved in time (with time step size $\tau = 10^{-5}$) by the convergent numerical algorithm \eqref{eq:BDF} from \cite{MCF} with numerical surgery, as proposed in Algorithm~\ref{alg:numerical MCF with surgery}, with curvature thresholds $H_2 = 100$ and $H_3 = 200$. 

The figure shows snapshots of the numerical solution $\Ga_h[\bfx^n]$ and the curvature $H_h^n$ at eight different time steps (see labels). A pinch singularity is resolved by numerical surgery at time $T_1 = 0.08242$ (third row), while the two remaining parts go extinct at $T_2 = 0.09967$ (bottom row). The evolution of the mean curvature $H_h^n$ can be observed using the colorbars.

\begin{figure}[htbp]
	\centering
	\includegraphics[height=4.5cm,clip,trim={0 0 0 0}]{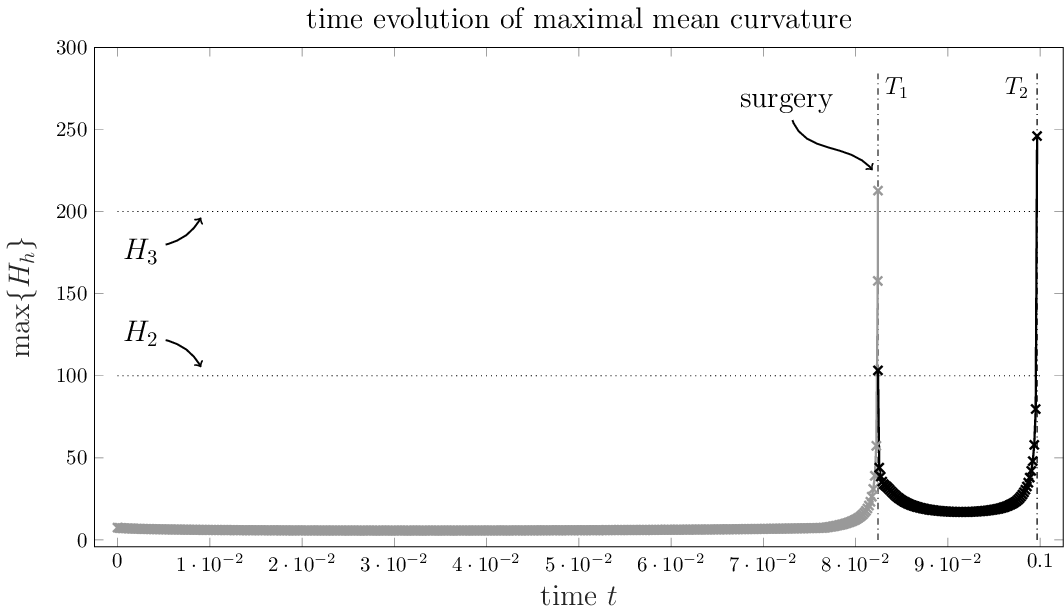}
	\includegraphics[height=4.55cm,clip,trim={0 0 0 0}]{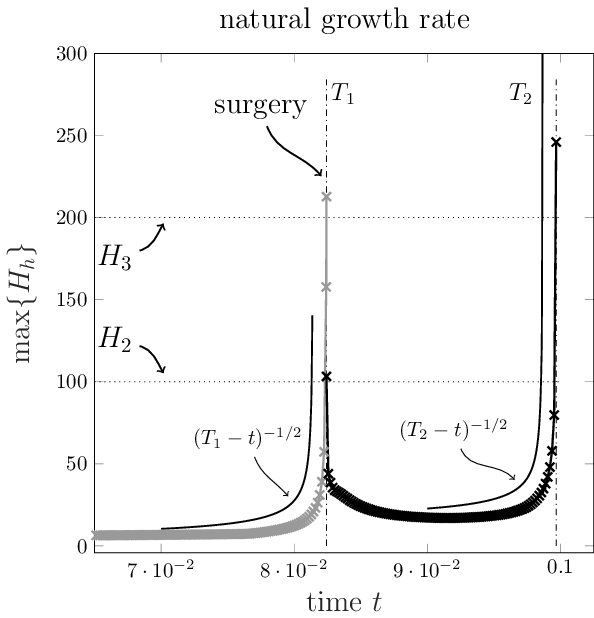}
	\caption{Time evolution of maximal mean curvature during numerical mean curvature flow with surgery for a dumbbell. Note the mean curvature drop due to surgery at time $T_1 = 0.08242$ (note the grey and black colouring). \bbk The plot on the right-hand side shows the natural growth rate, cf.~\eqref{eq:natural growth rate}. \ebk }
	\label{fig:mean_curvature}
\end{figure}

We further report on the time evolution of the maximal mean curvature. Figure~\ref{fig:mean_curvature} plots $\max_{\Ga_h[\bfx^n]}\{H_h^n\}$ against time over the interval $[0,T_2] = [0,0.09967]$. The curvature thresholds $H_2$ and $H_3$, and the times $T_1$ and $T_2$ are marked. We plotted the curves coloured grey and black before and after surgery, respectively, in order to better distinguish pre- and post-surgery phases. Due to the discrete nature of the flow we do not have $H_h^n \leq H_3$ throughout the flow. This could easily be overcome by performing surgery for the surface of the previous time step.

\begin{figure}[htbp]
	\centering
	\includegraphics[width=0.32\textwidth,trim={377 40 382 50}]
	{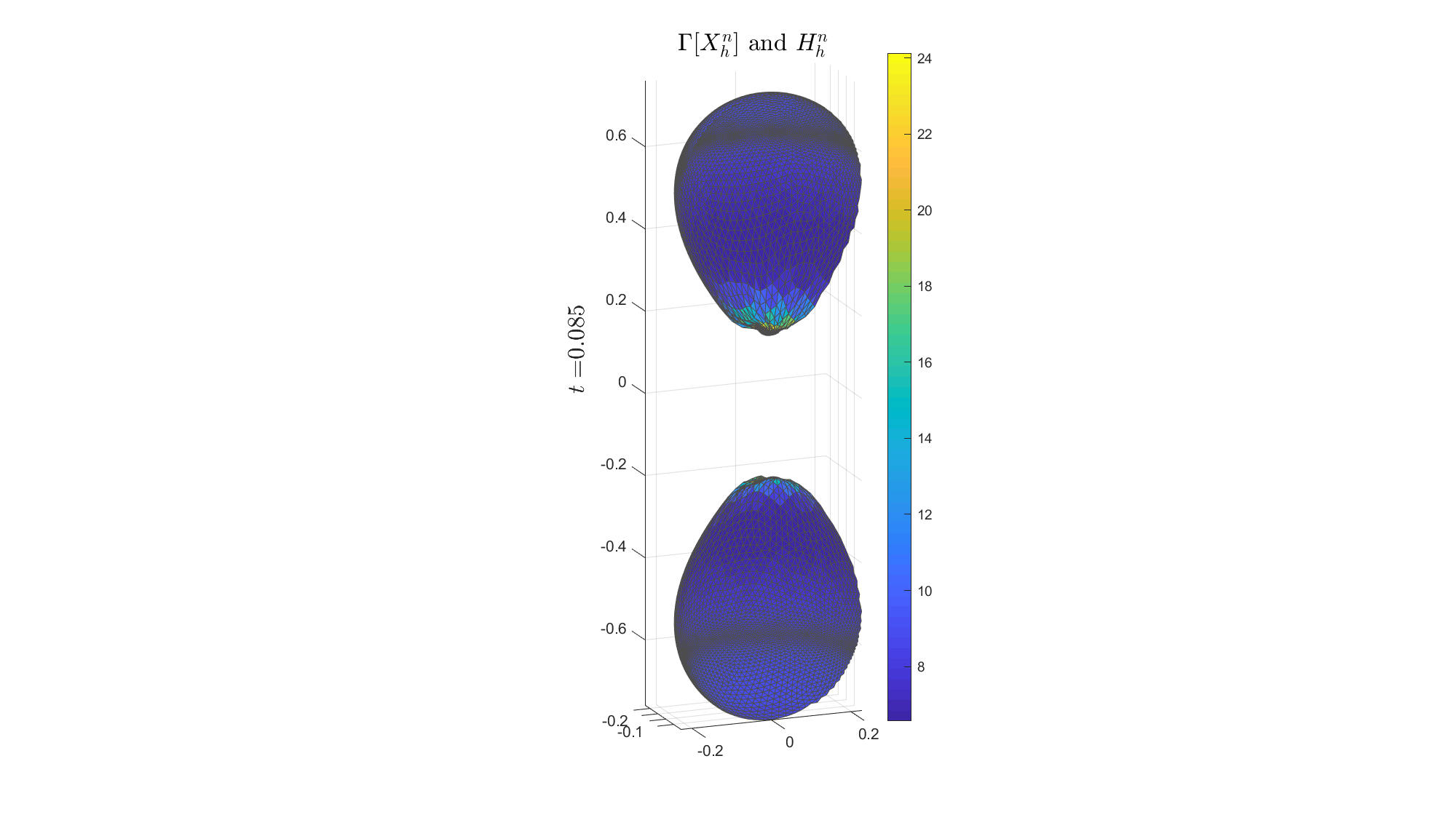}
	\includegraphics[width=0.32\textwidth,clip,trim={370 40 390 20}]
	{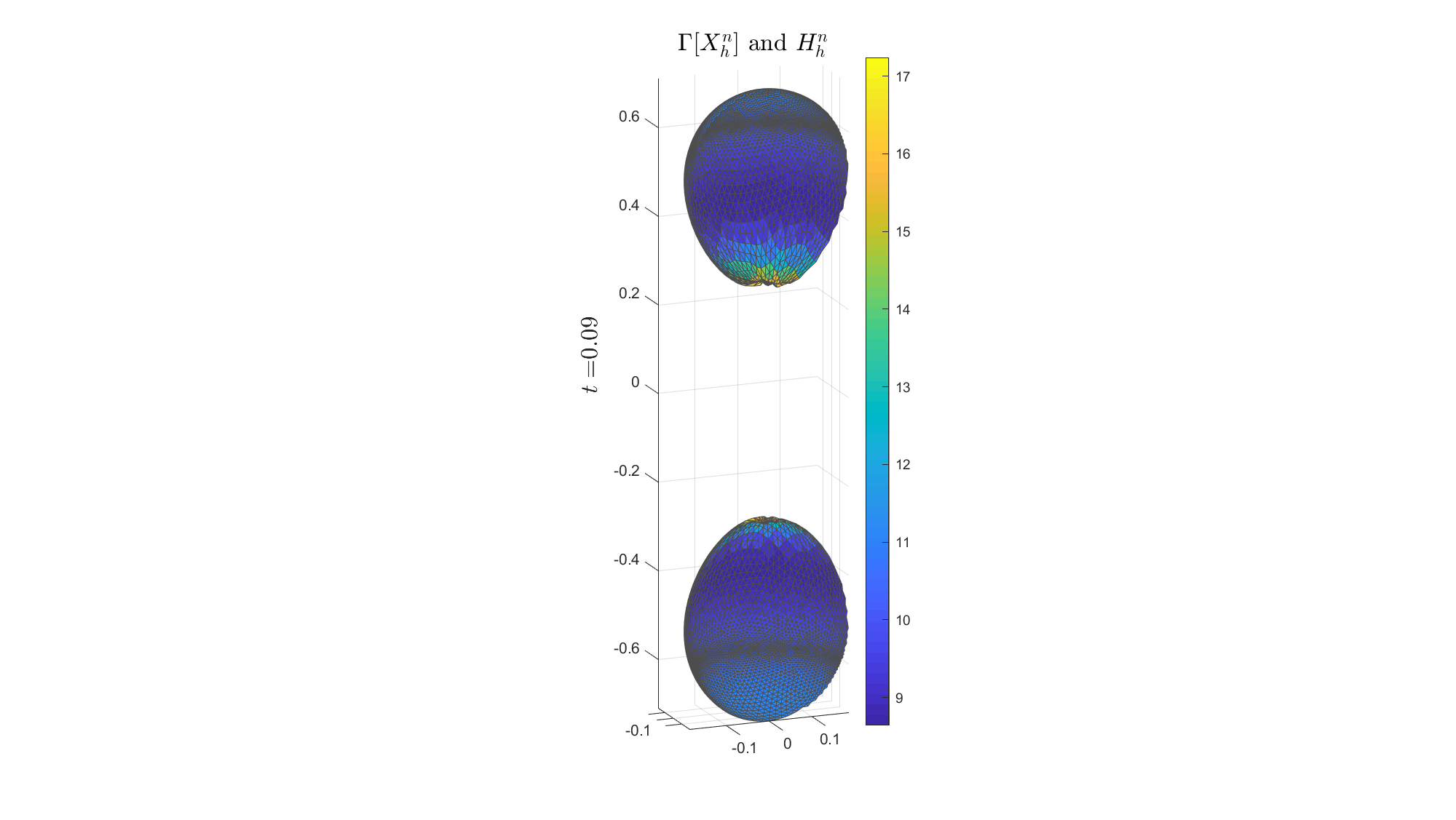}
	\includegraphics[width=0.32\textwidth,clip,trim={370 40 391 20}]
	{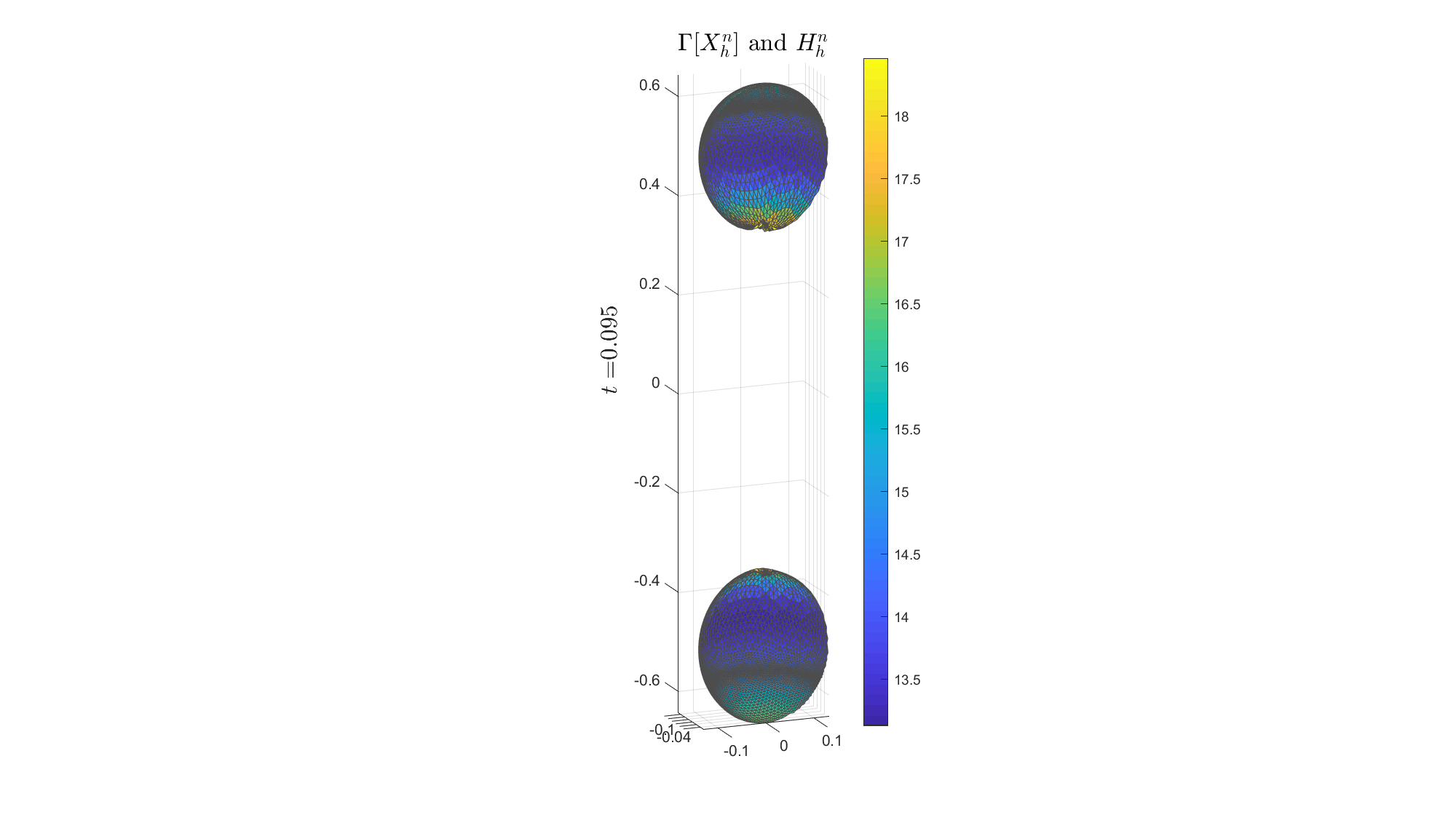}
	
	\includegraphics[width=0.32\textwidth,clip,trim={370 40 412 10}]
	{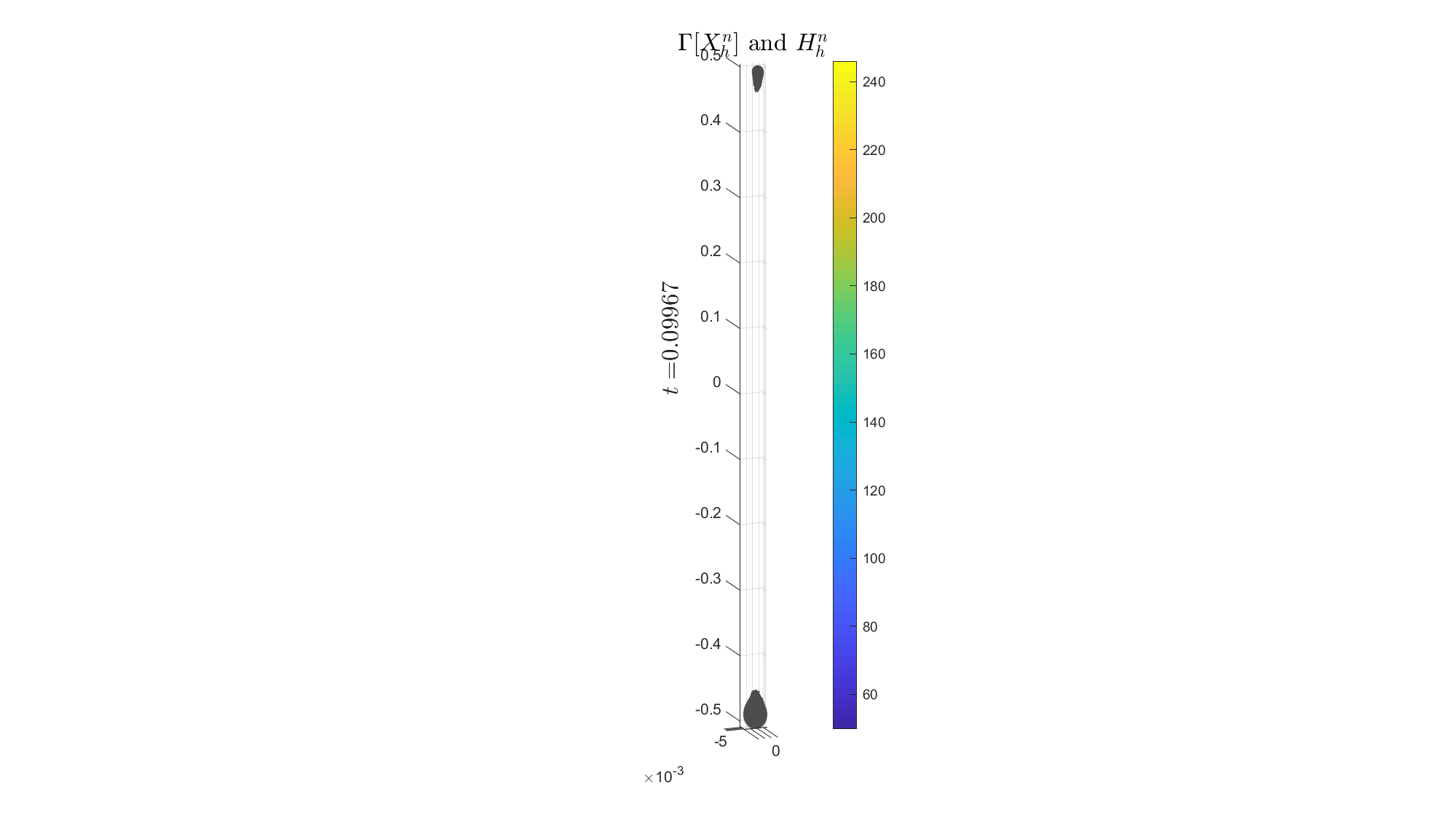}
	\includegraphics[width=0.32\textwidth,clip,trim={350 0 345 20}]
	{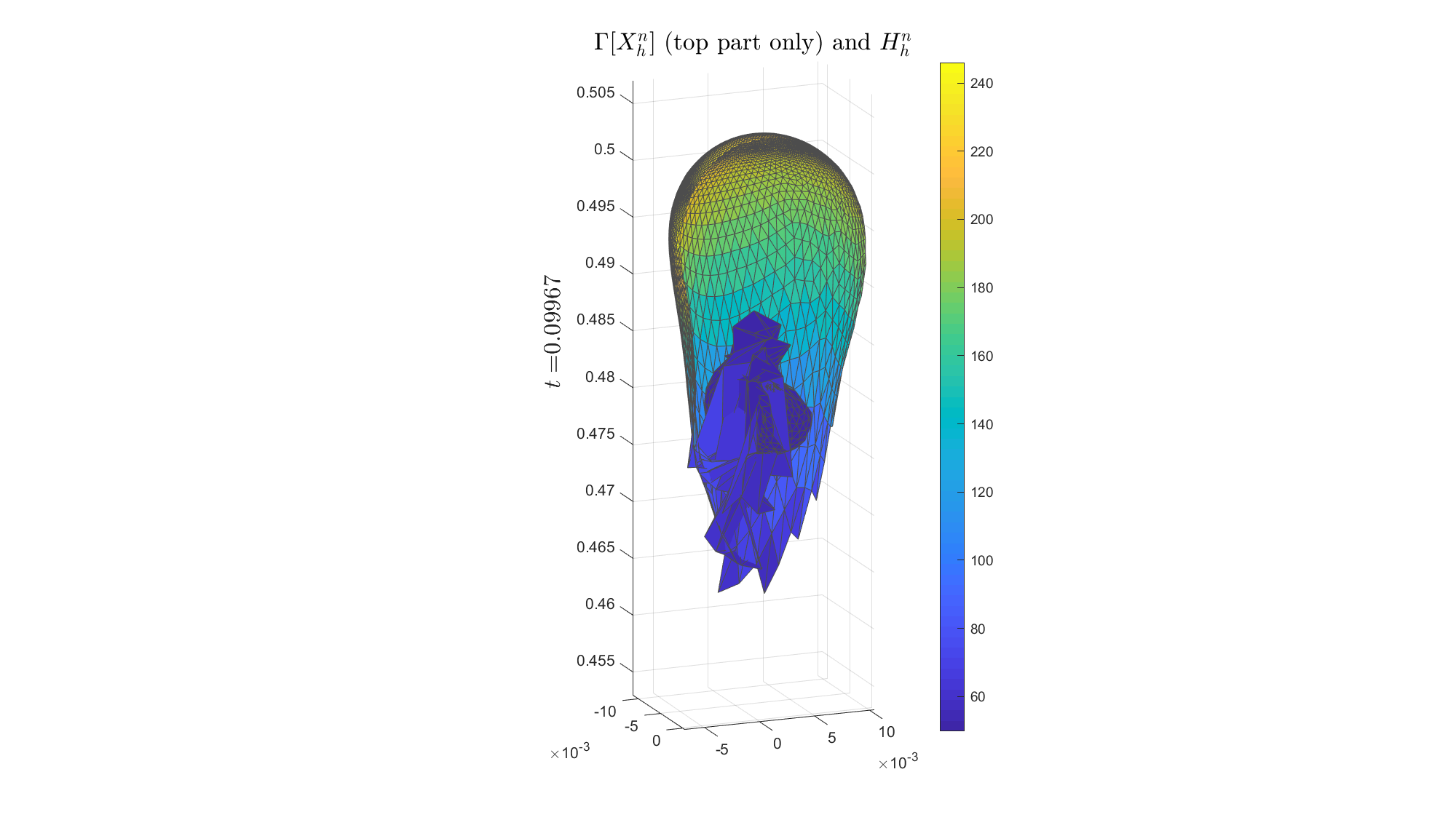}
	\caption{The cross sections show that the two spherical parts indeed shrink to a \emph{round point} at $T_2 = 0.09967$. Odd mesh behaviour is only observed at the singularity (see final image; note the scaling).}
	\label{fig:cross_sections}
\end{figure}

In Figure~\ref{fig:cross_sections} we report on cross sections of the numerical solution, in order to demonstrate that spikes or other odd mesh artefacts are not developed inside the surfaces before the extinction time (last row). The right-hand side image in the last row shows a zoom-in for the top part (note the axes are of magnitude $10^{-3}$). 

\pagebreak
\subsubsection*{The numerical surgery}

\begin{figure}[htbp]
	\centering
	\includegraphics[width=1\textwidth,clip,trim={110 30 90 20}]
	{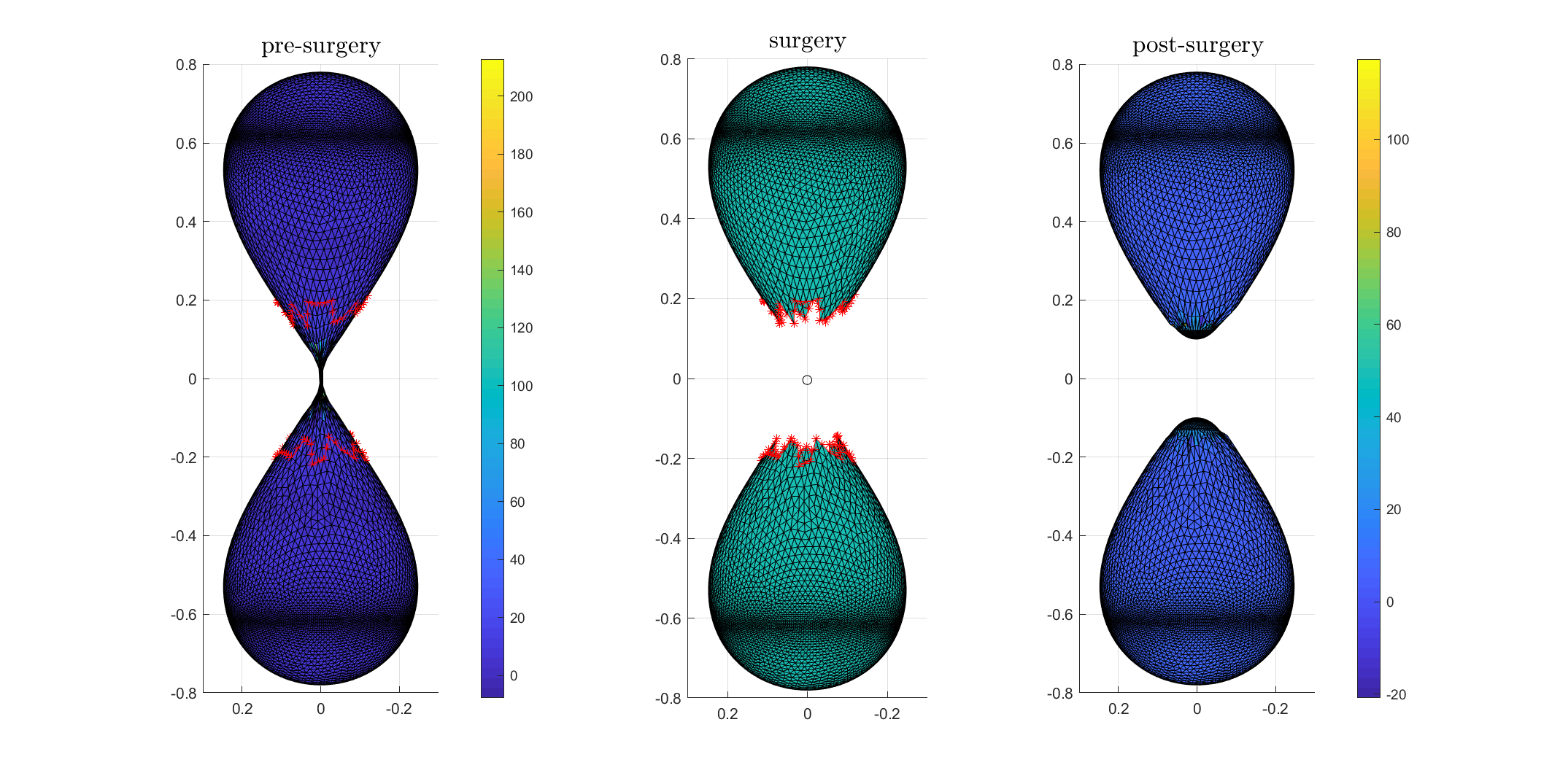}
	\includegraphics[width=1\textwidth,clip,trim={100 100 80 100}]
	{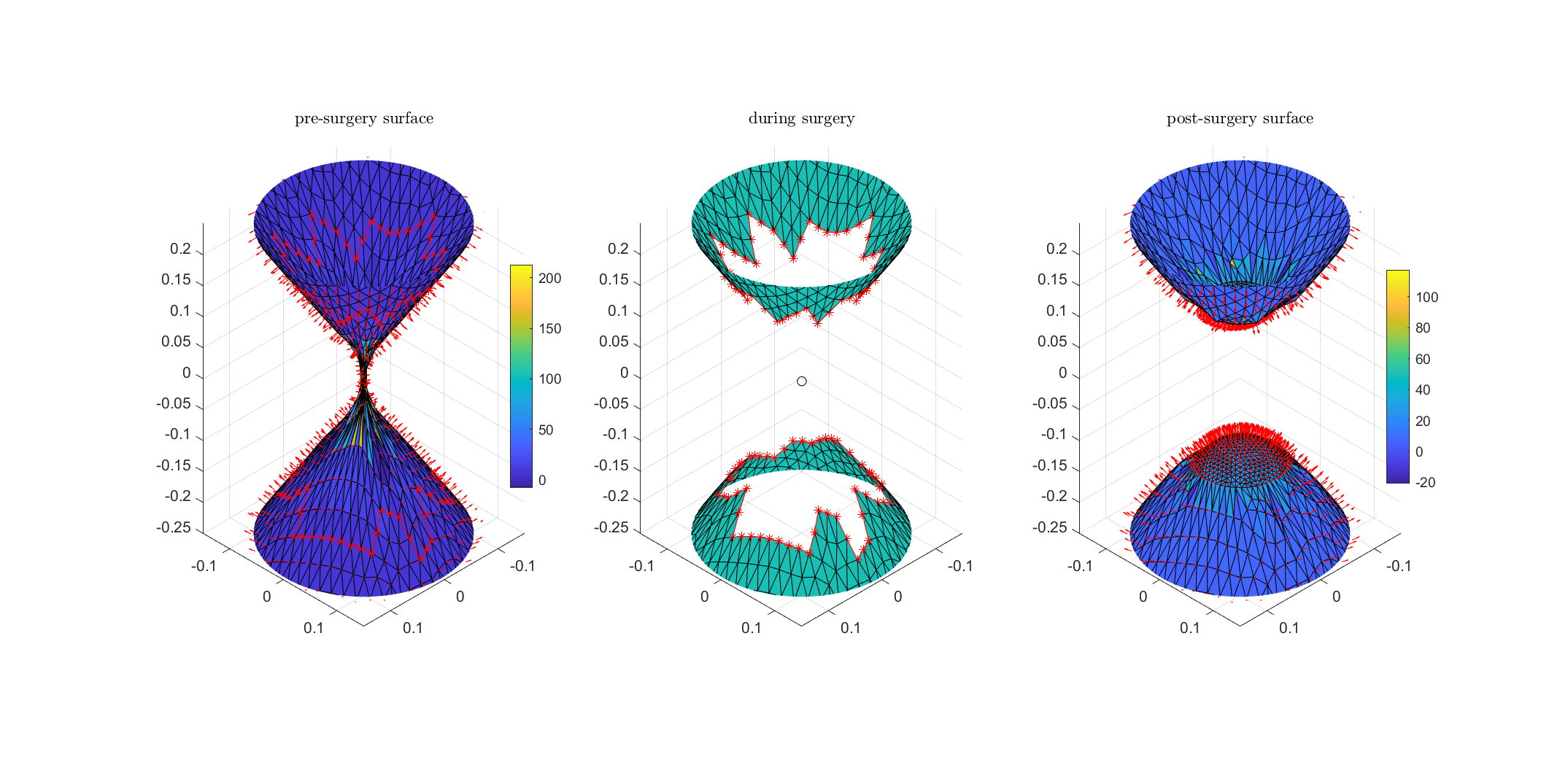}
	\caption{The steps of numerical surgery Algorithm~\ref{alg:numerical surgery} for the dumbbell example, when $\max\{H_h^n\}$ exceeds $H_3$ at time $T_1 = 0.08242$, reporting on $H_h$ and $\nu_h$ as well.}
	\label{fig:numerical_surgery}
\end{figure}

In Figure~\ref{fig:numerical_surgery} we report on the numerical surgery, i.e.~on Algorithm~\ref{alg:numerical surgery}, occurring at time $T_1 = 0.08242$.
The figure shows (close-up view in the bottom row) the discrete surface $\Ga_h[\bfx^n]$ and the discrete mean curvature $H_h^n$ (see colorbars). The two columns on the left-hand side show the surgery boundary (red stars), while the bottom row also shows the normal vector field before and after the surgery (arrows in red).

\subsubsection{A comparison with the level set method}

We reproduced Figure~\ref{fig:MCF_with_surgery} using the level set method of \cite{OsherSethian_1988}---which automatically resolves singularities, using the \texttt{ToolboxLS} of \cite{Mitchell_1,Mitchell_2} (using quadratic approximations in space and time\footnote{See \texttt{upwindFirstENO2} and \texttt{odeCFL2} in \cite[Section~3.4.1 \& 3.5.1]{Mitchell_2}.}). 

\begin{figure}[htbp]
	\centering
	\includegraphics[width=\textwidth,clip,trim={0 45 0 10}]
	{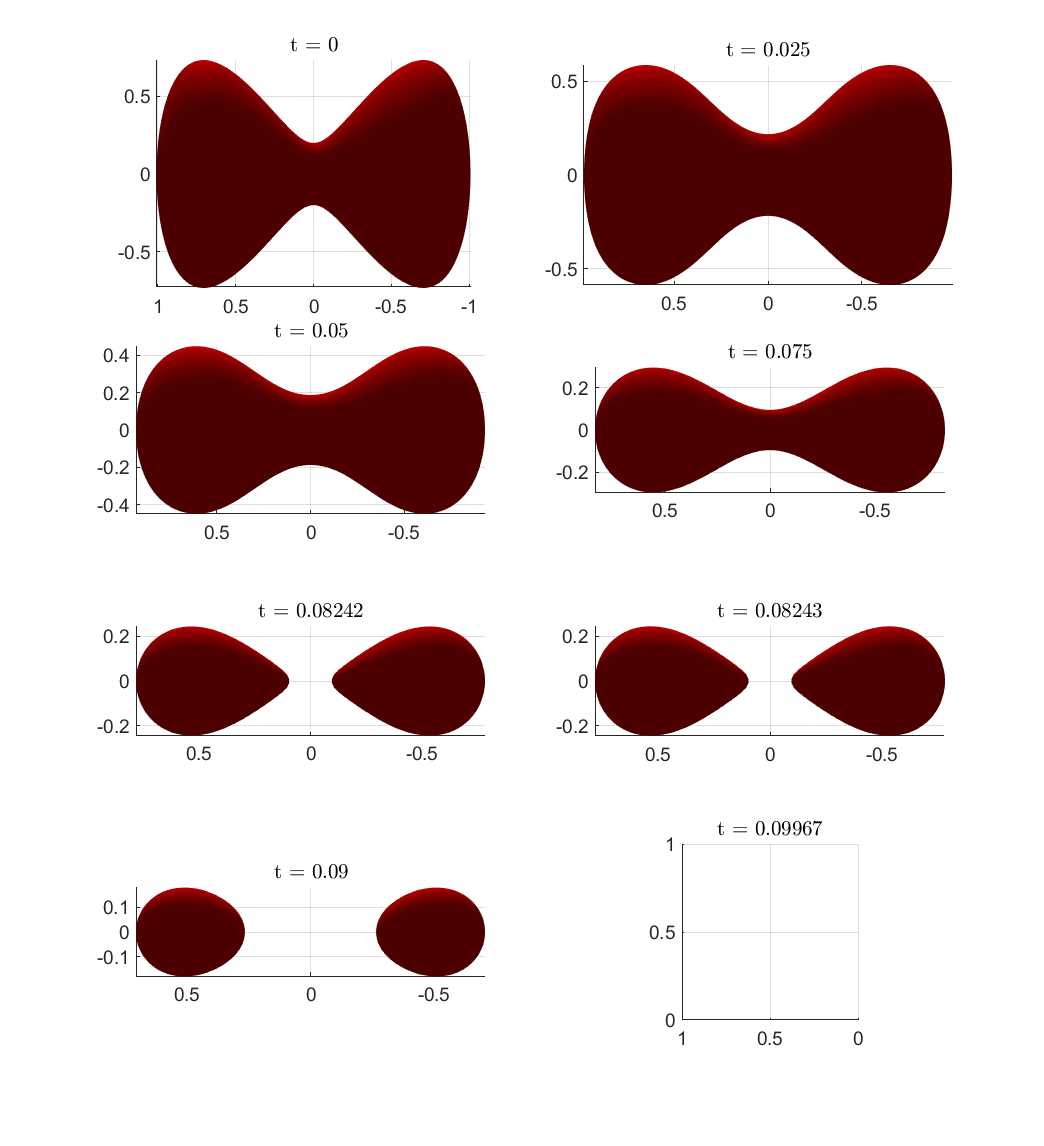}
	\caption{The level set method version of Figure~\ref{fig:MCF_with_surgery} at identical time steps. The plots were generated using \texttt{ToolboxLS} \cite{Mitchell_1,Mitchell_2}.}
	\label{fig:Dumbbell_singular_ToolboxLS}
\end{figure}

The obtained surface approximations are shown in Figure~\ref{fig:Dumbbell_singular_ToolboxLS}. For easy comparison the plots show the numerical solutions at the exact same time-steps as Figure~\ref{fig:MCF_with_surgery}. 

The two numerical solutions are extremely close to each other. 
Note that the level set method resolves the pinch singularity earlier than the Algorithm~\ref{alg:numerical MCF with surgery}, and the remaining components also vanish earlier. Note that, to our knowledge, no error estimates are known for the level set method.

\subsection{The pinch-off of a sphere from a torus}

For this experiment the initial surface, a sphere and a torus connected by a thin neck, is given as a constant distance product surface. 
Let the functions $d_{\text{torus}}$, and $d_{\text{S}_1}$, $d_{\text{S}_2}$ be, respectively, the distance functions of a torus with radii $R=4$ and $r=2$, and of spheres with centres $(4,0,2.5)$, $(4,0,5.25)$ and radius $0.5$, $2.5$. Then $\Ga^0$ is given as the zero level set of the function
\begin{equation}
\label{eq:torus_and_sphere}
d(x) = d_{\text{torus}}(x) \cdot d_{\text{S}_1}(x) \cdot d_{\text{S}_2}(x) - 0.07 .
\end{equation} 

\begin{figure}[htbp]
		\centering
		\includegraphics[height=0.18\textheight,clip,trim={82 10 48 8}]
		{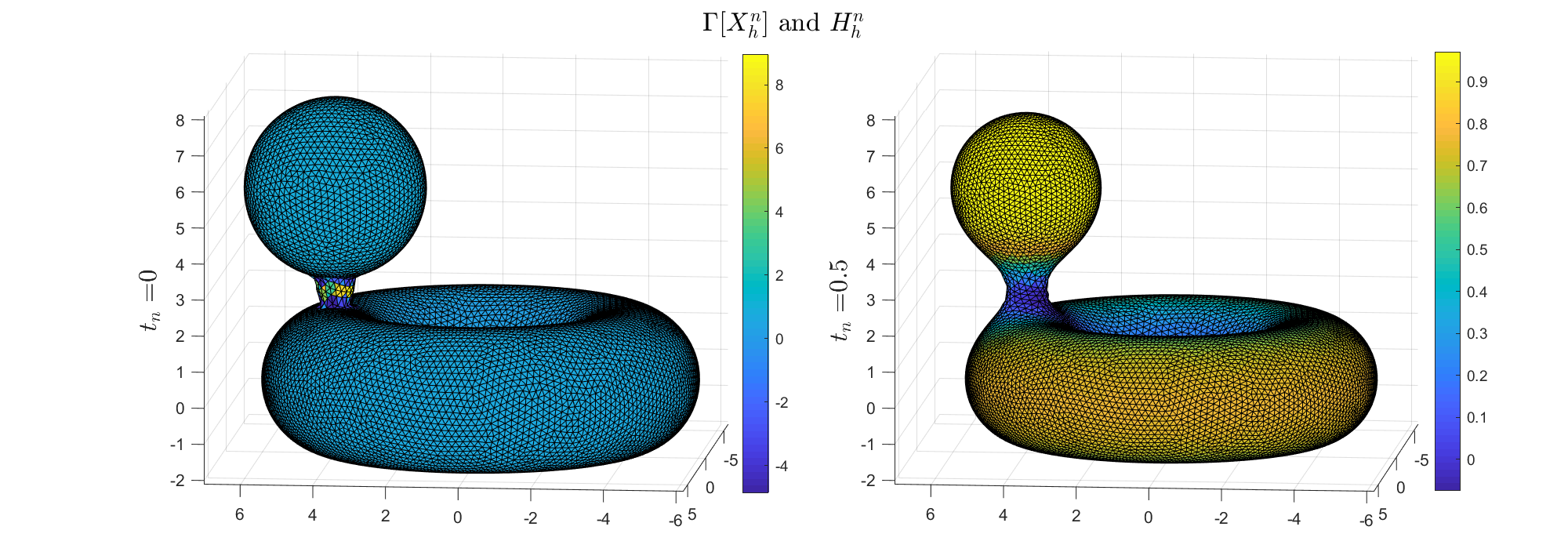}
		\includegraphics[height=0.18\textheight,clip,trim={82 18 48 25}]
		{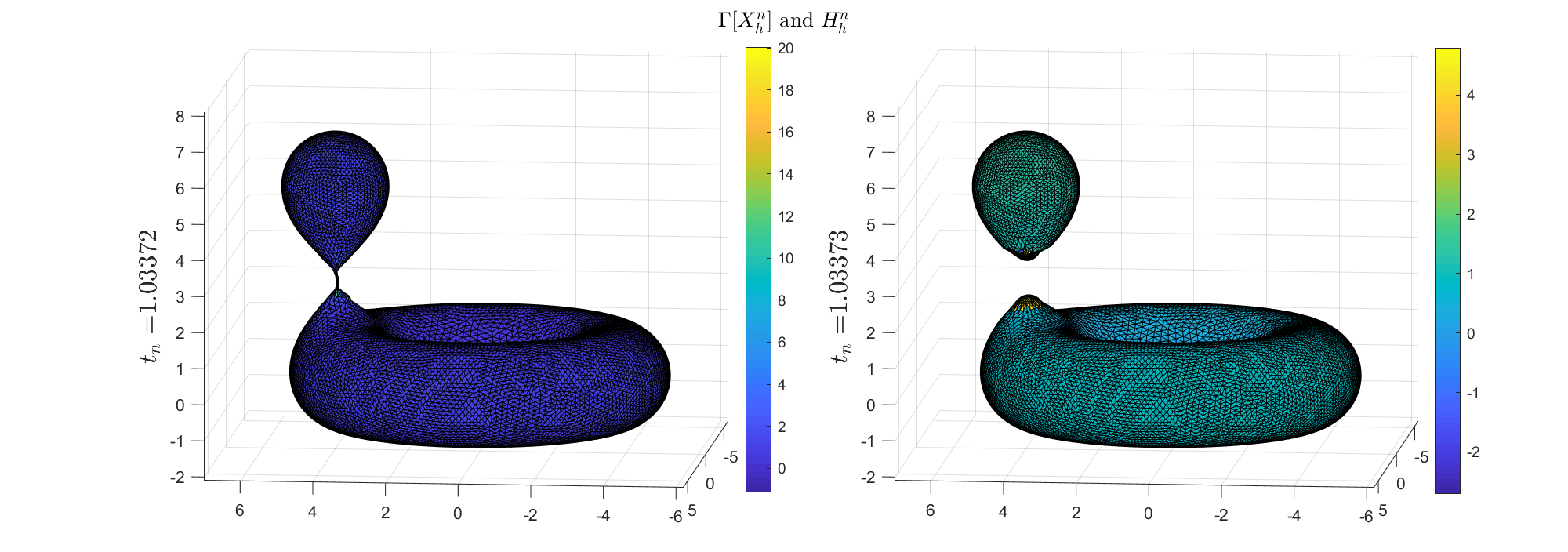}
		\includegraphics[height=0.18\textheight,clip,trim={82 18 48 25}]
		{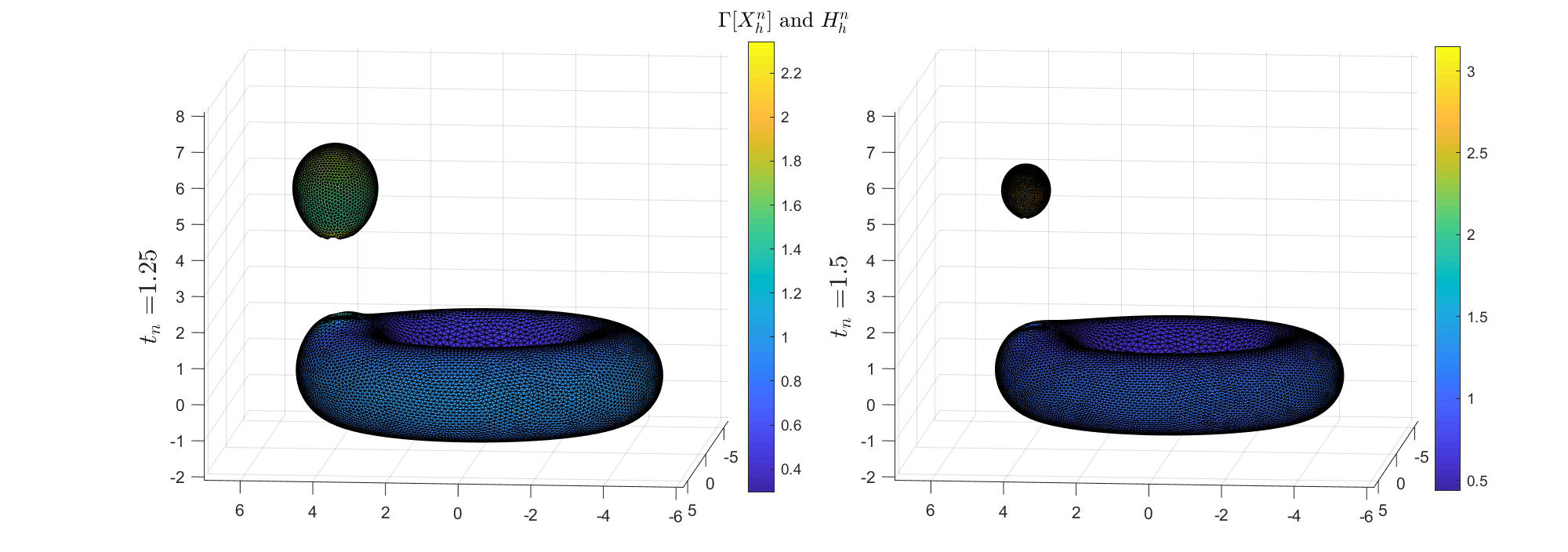}
		\includegraphics[height=0.18\textheight,clip,trim={82 18 48 25}]
		{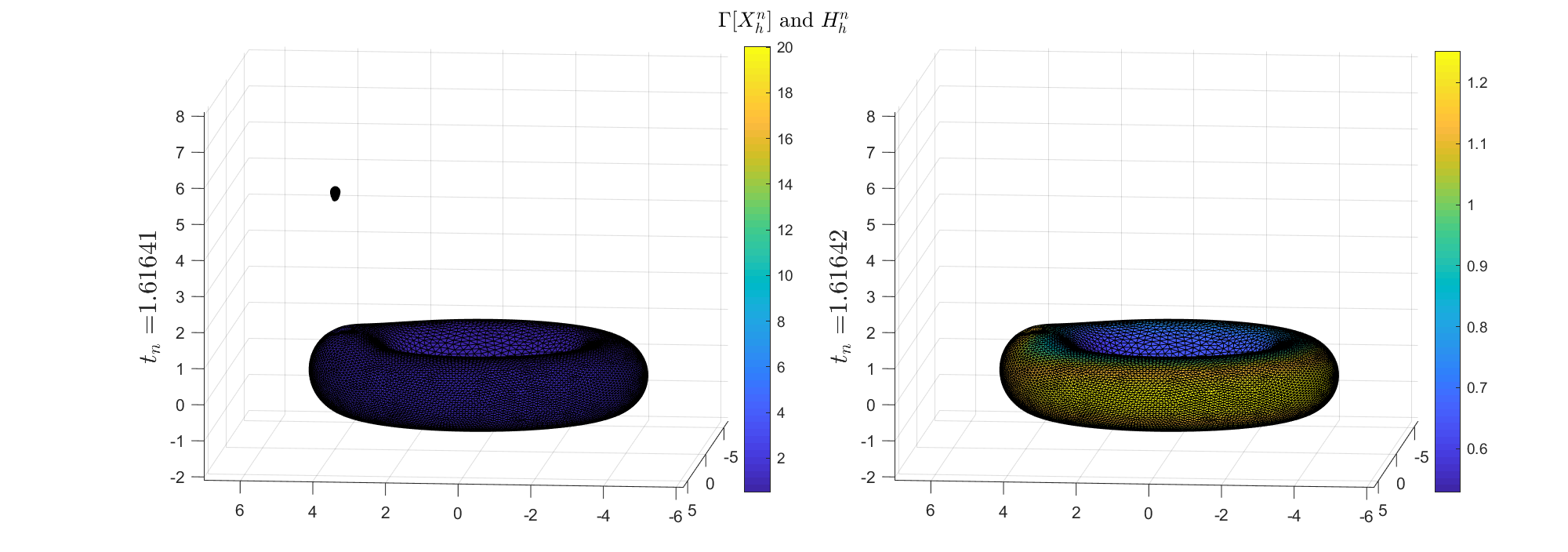}
		\includegraphics[height=0.18\textheight,clip,trim={82 18 48 25}]
		{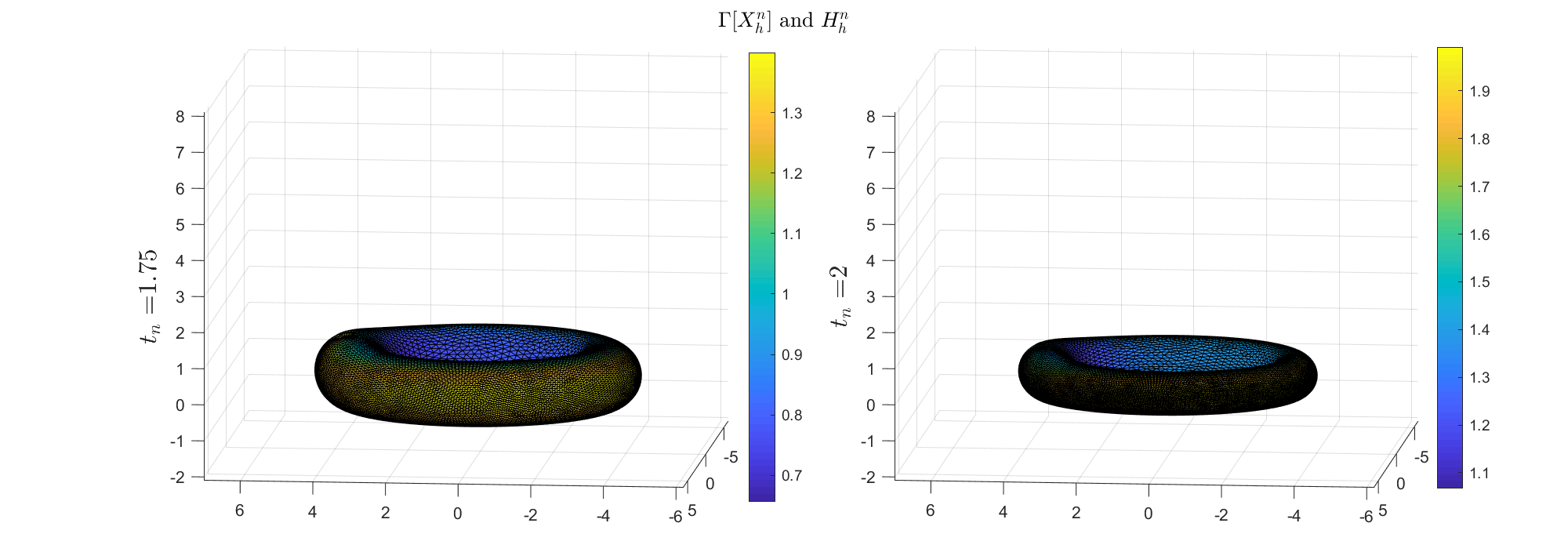}
	\caption{The numerical solution, $\Ga[X_h^n]$, and $H_h^n$ at various time steps $t_n$, of mean curvature flow with surgery of a torus and sphere joined by a think neck.}
	\label{fig:Torus_and_Sphere}
\end{figure}

\pagebreak
Figure~\ref{fig:Torus_and_Sphere} reports on numerical mean curvature flow with surgery for this torus--sphere initial surface (given by \eqref{eq:torus_and_sphere}).\footnote{A video of this experiment is available at \href{https://kovacs.app.uni-regensburg.de/research.html}{\texttt{https://kovacs.app.uni-regensburg.de/research.html}}.} The interpolation of the initial surface $\Ga^0$ (see top left in Figure~\ref{fig:Torus_and_Sphere}) is evolved in time by Algorithm~\ref{alg:numerical MCF with surgery} ($\tau = 10^{-5}$, and $\text{dof} = 11496$ and $h = 0.30815$ initially). The curvature thresholds are $H_2 = 15$, and $H_3 = 20$. Recall with a given maximal mesh width $h$ one cannot meaningfully represent a discrete surface with arbitrarily high curvature. Note that the maximal mesh width of initial surface $\Ga_h^0$ is roughly $10$ times larger than for the dumbbell example. 

Figure~\ref{fig:Torus_and_Sphere} shows snapshots of the numerical solution $\Ga_h[\bfx^n]$ and the curvature $H_h^n$ at different time steps (see labels): a pinch singularity is resolved by numerical surgery at time $T_1 = 1.03372$ (2nd row), then the spherical part shrinks to a round point at $T_2 = 1.61641$ and is hence removed (4th row). The evolution of the mean curvature $H_h^n$ can be observed using the colorbars.

\begin{figure}[htbp]
	\centering
	\includegraphics[width=1\textwidth,clip,trim={0 0 0 0}]{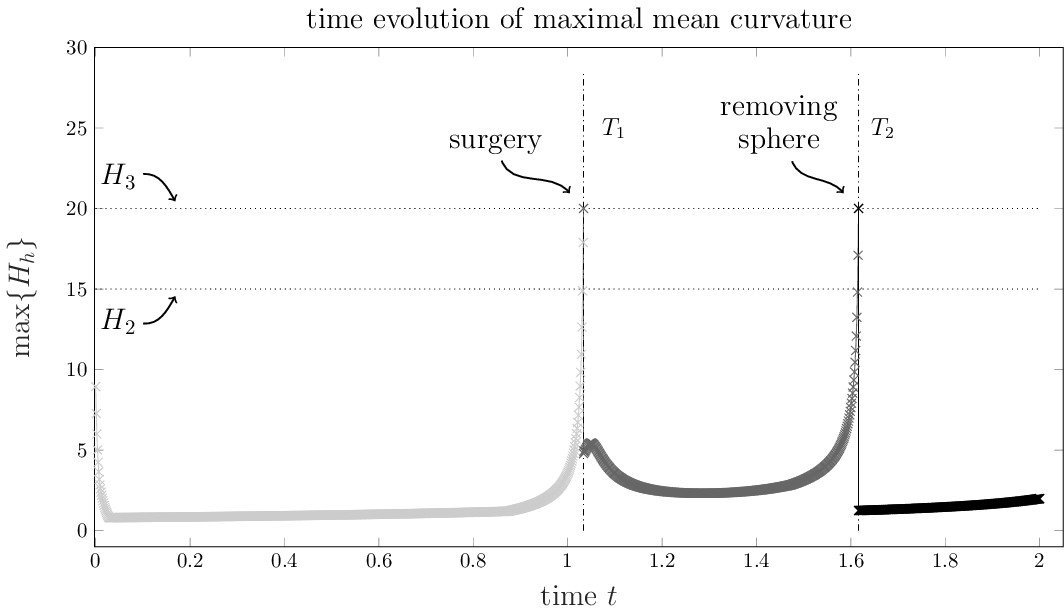}
	\caption{Time evolution of maximal mean curvature during numerical mean curvature flow with surgery of a torus and sphere joined by a think neck. (Note the grey-scaled colouring separating the flow between surgeries.)}
	\label{fig:mean_curvature_TS}
\end{figure}

We again report on the time evolution of the maximal mean curvature, Figure~\ref{fig:mean_curvature_TS}. The curvature thresholds $H_i$, and the surgery times $T_j$ are marked again.

\begin{figure}[htbp]
	\centering
	\includegraphics[width=1\textwidth,clip,trim={110 145 60 130}]
	{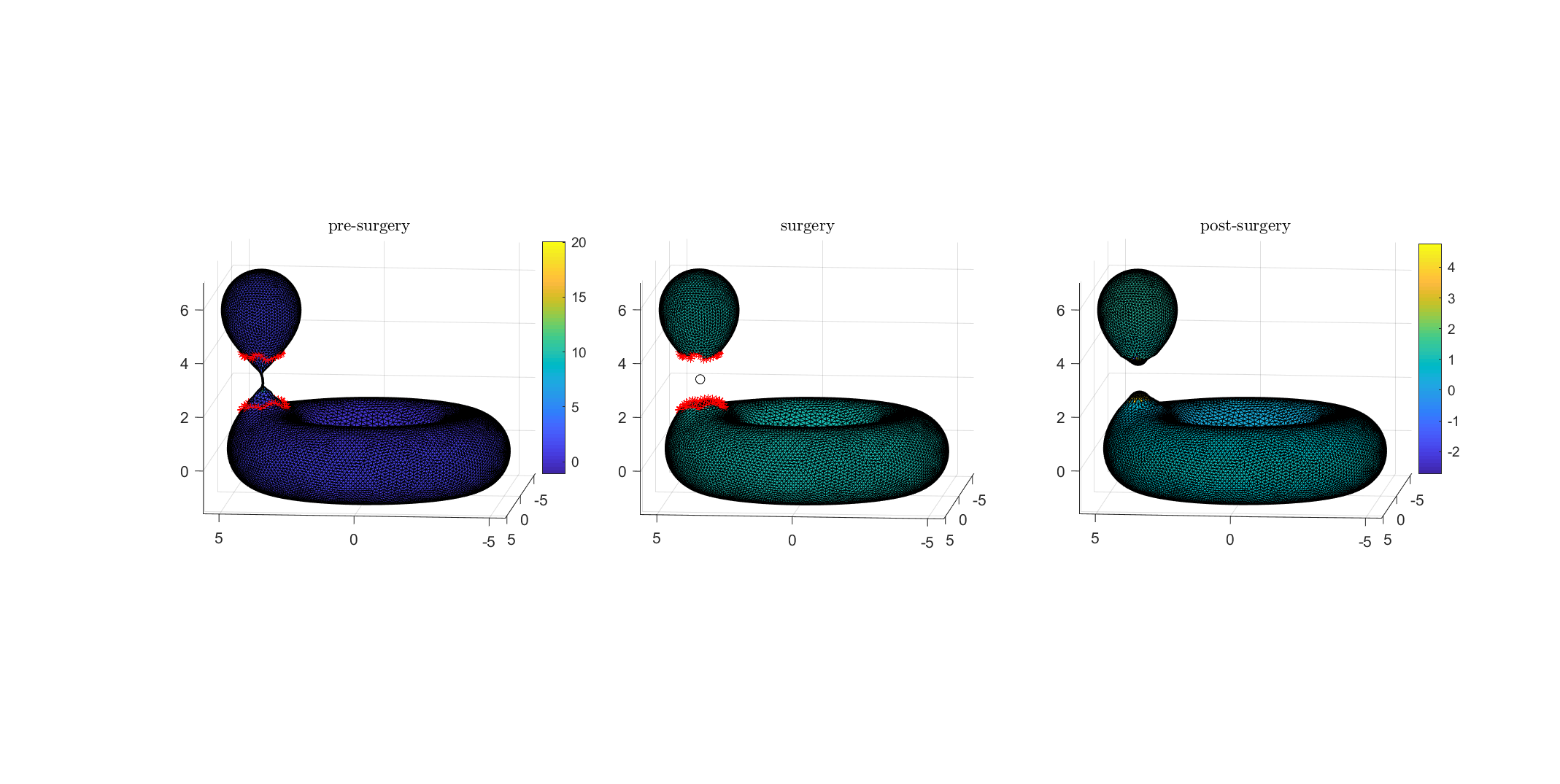}
	\includegraphics[width=1\textwidth,clip,trim={100 100 60 100}]
	{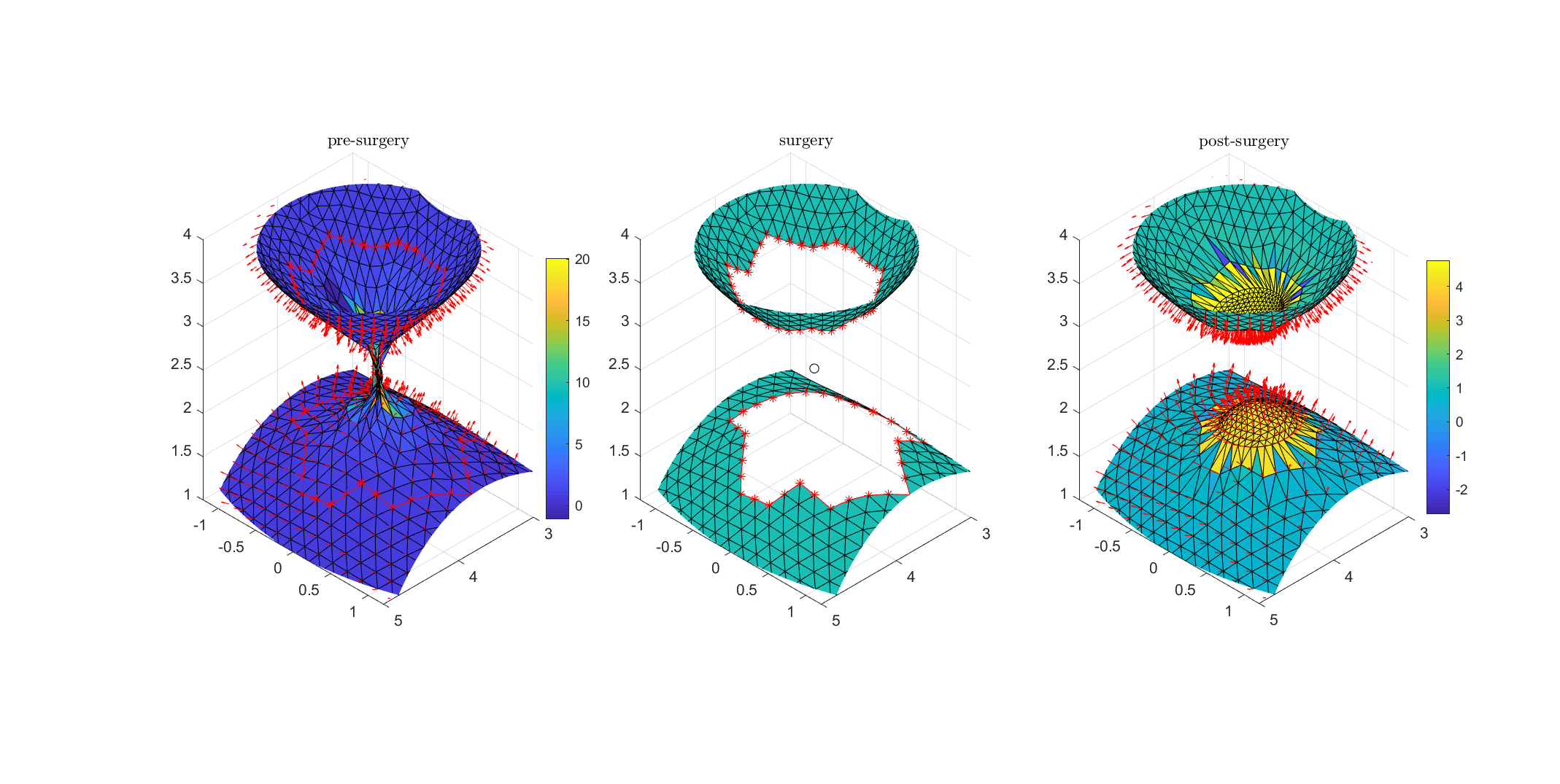}
	\includegraphics[width=1\textwidth,clip,trim={110 145 60 130}]
	{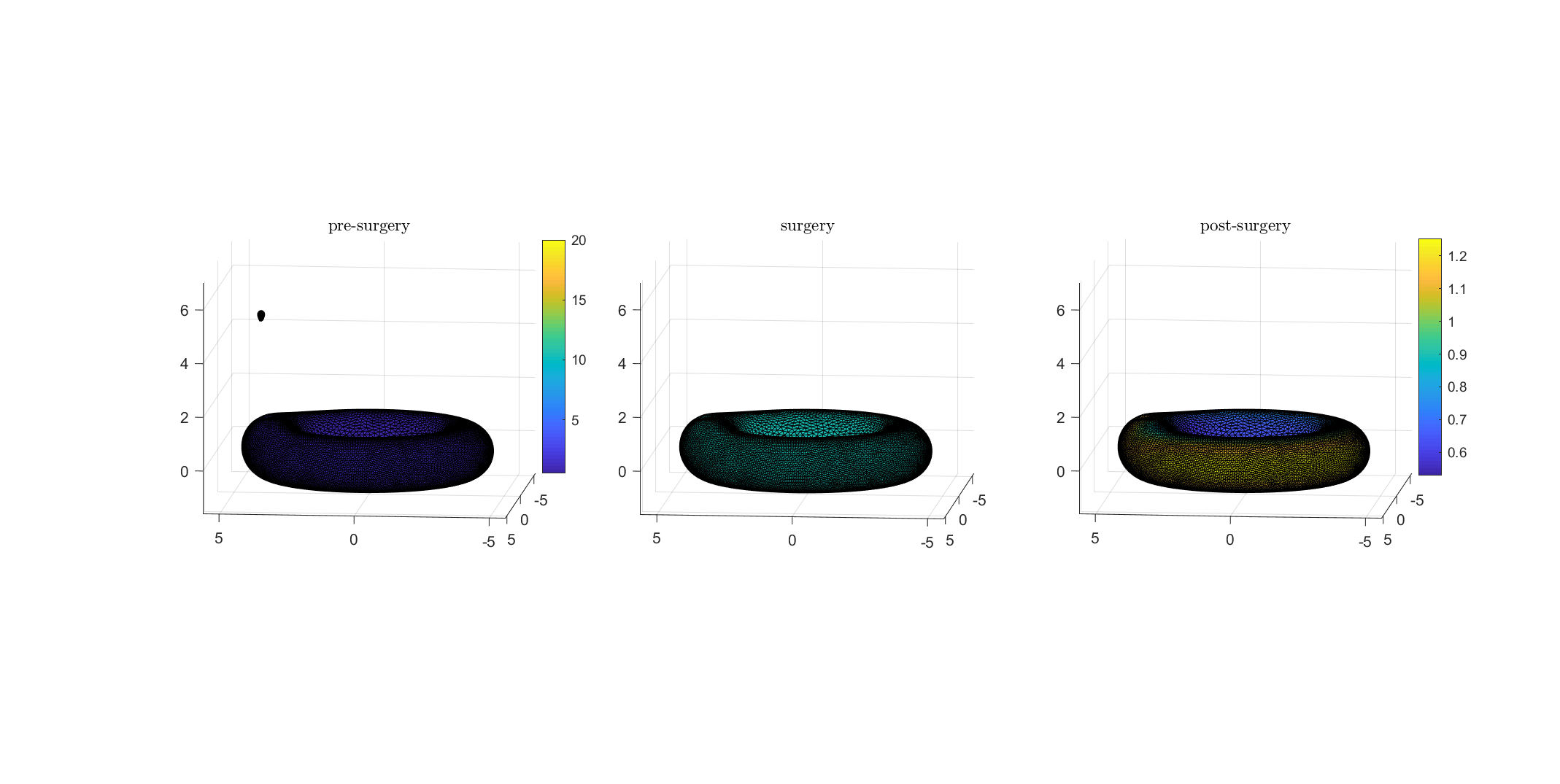}
	\caption{The steps of numerical surgery Algorithm~\ref{alg:numerical surgery} for the torus--sphere example, when $\max\{H_h^n\}$ exceeds $H_3$ at time $T_1 = 1.03372$ (first two rows), and the removal of a small spherical part when $\max\{H_h^n\}$ exceeds $H_3$ again at time $T_2 = 1.61641$ (bottom row).}
	\label{fig:numerical_surgery_TS}
\end{figure}

Similarly as before, in Figure~\ref{fig:numerical_surgery_TS} we report on the numerical surgeries occurring at time $T_1 = 1.03372$ and $T_2 = 1.61641$ for the joined torus and sphere example. At $T_1$ a pinch singularity is resolved much like before (first two rows), while at $T_2$ the pinched-off top portion have shrinked to spherical round point, and hence was removed (bottom row).

\pagebreak

\subsection{Torus merging into a sphere}
\label{section:fat torus}

We implemented the generalised numerical surgery for \emph{non-mean convex} surfaces, see Remark~\ref{remark:fat torus}. 

For this experiment the initial surface is a thick torus---which will merge into a sphere---is given as $\Ga^0 = \{x \in \R^3 \mid d(x) = 0 \}$ with the signed distance function
\begin{equation}
\label{eq:distance - fat torus}
	d(x) = \big( (x_1^2 + x_2^2)^{1/2} - R \big)^2 + x_3^2 - r^2 , \qquad \text{with} \quad R = 3 \quad \text{and} \quad r = 2.75.
\end{equation}

Similarly as before, we report on the evolution of the surface the mean curvature, and also on the surgery itself.

\begin{figure}[htbp]
	\centering
	\includegraphics[height=0.18\textheight,clip,trim={82 93 53 00}]
	{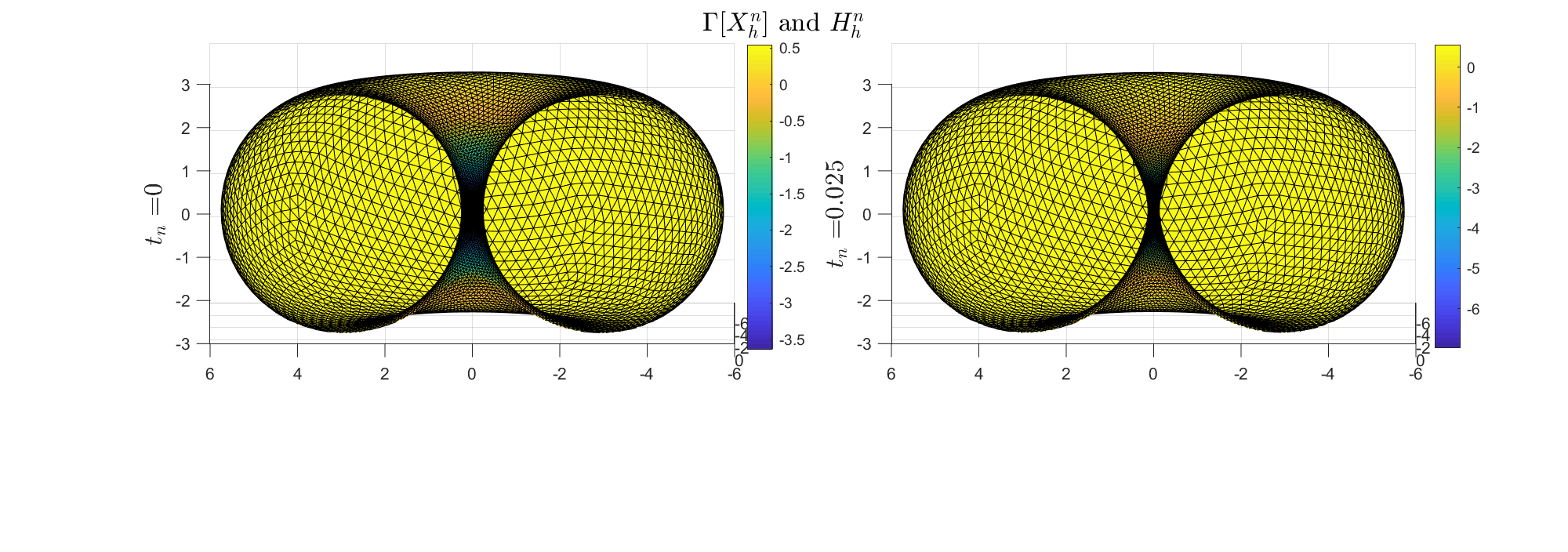}
	\includegraphics[height=0.18\textheight,clip,trim={80 50 0 50}]
	{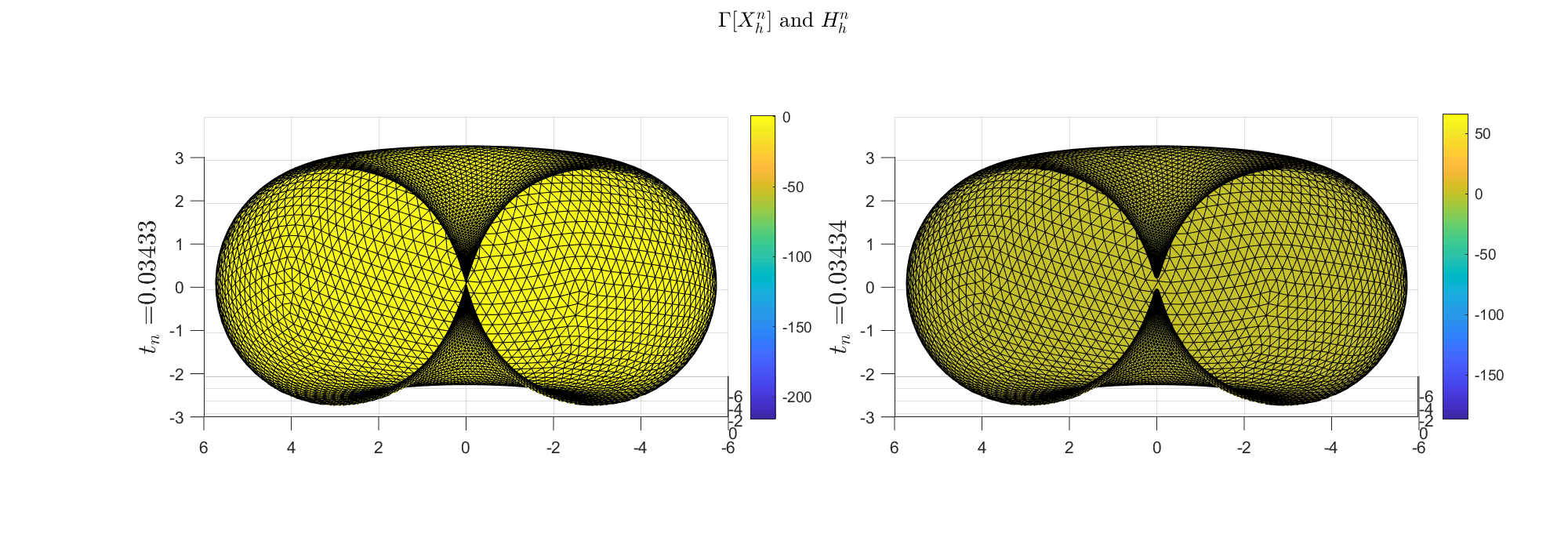}
	\includegraphics[height=0.18\textheight,clip,trim={80 50 0 50}]
	{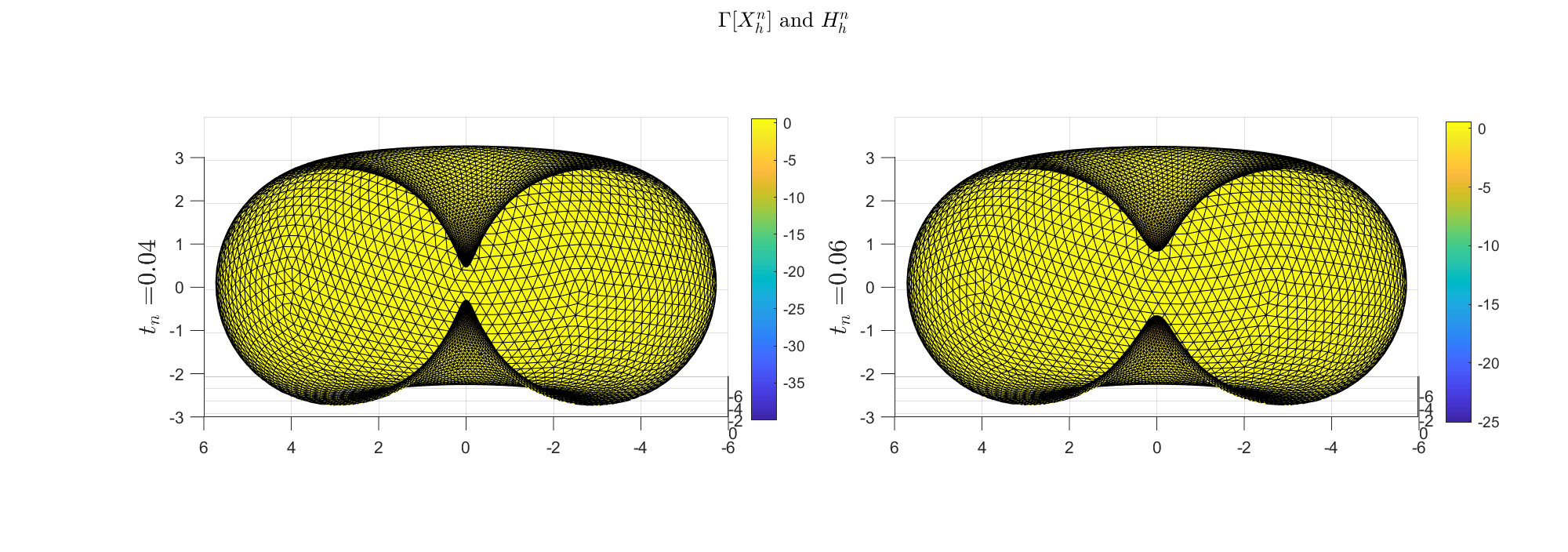}
	\includegraphics[height=0.18\textheight,clip,trim={80 50 0 50}]
	{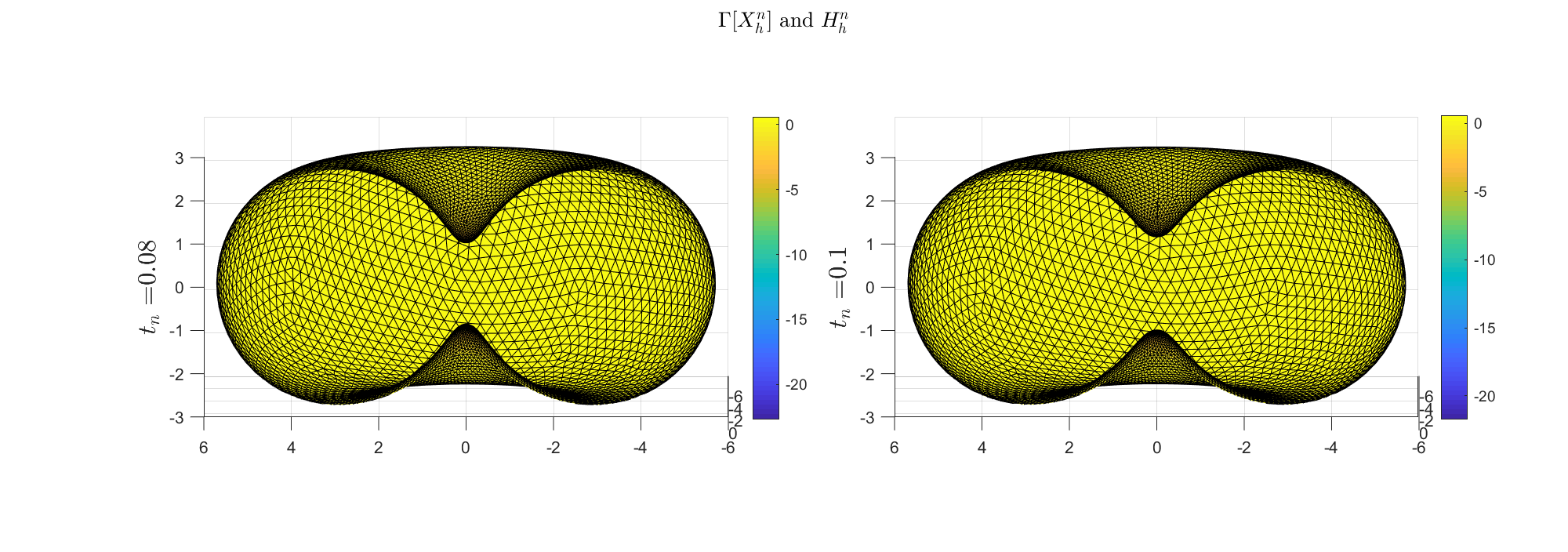}
	\caption{The numerical solution, $\Ga[X_h^n]$, and $H_h^n$ at various time steps $t_n$, of mean curvature flow with surgery of a torus merging into a sphere.}
	\label{fig:Torus_fatG}
\end{figure}

Figure~\ref{fig:Torus_fatG} reports on numerical mean curvature flow with surgery for the torus \eqref{eq:distance - fat torus} as initial surface.\footnote{A video of this experiment is available at \href{https://kovacs.app.uni-regensburg.de/research.html}{\texttt{https://kovacs.app.uni-regensburg.de/research.html}}.} 
To properly resolve the thin hole of the torus, we used a graded mesh interpolating the initial surface $\Ga^0$ (see top left in Figure~\ref{fig:Torus_fatG}). The initial surface $\Ga_h^0$ is evolved in time by Algorithm~\ref{alg:numerical MCF with surgery} ($\tau = 10^{-5}$, and $\text{dof} = 11488$ with $h_{\max} = 0.7145$ and $h_{\min} = 0.0326$ initially). The curvature thresholds are $H_2 = 100$ and $H_3 = 200$.

Figure~\ref{fig:Torus_fatG} shows snapshots of the numerical solution $\Ga_h[\bfx^n]$ and the curvature $H_h^n$ at different time steps (see labels): the singularity when the torus merges into a sphere is resolved by numerical surgery at time $T_1 = 0.03433$ (2nd row), then the remaining spikes start to collapse and the spherical surface starts to shrink. The evolution of the mean curvature $H_h^n$ can be observed using the colorbars.

\begin{figure}[htbp]
	\centering
	\includegraphics[width=1\textwidth,clip,trim={0 0 0 0}]{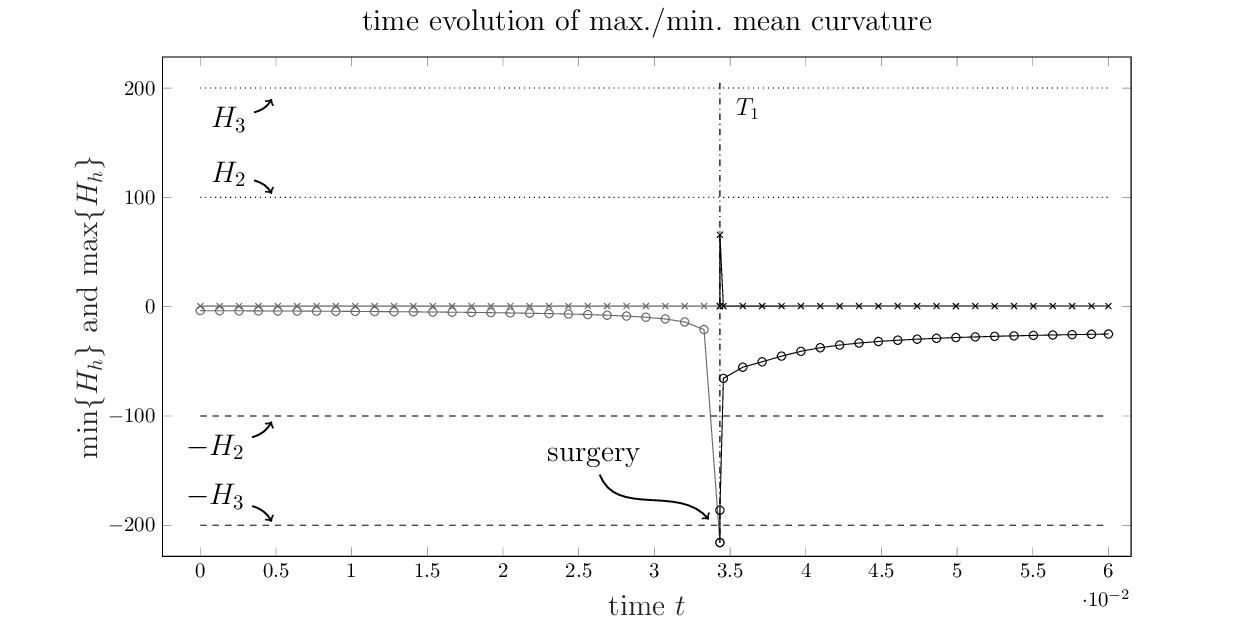}
	\caption{Time evolution of $\min\{ H_h^n \}$ and $\max\{ H_h^n \}$ during numerical mean curvature flow with surgery of a torus merging into a sphere. (Note the grey-scaled colouring separating the flow between surgeries.)}
	\label{fig:mean_curvature_Torus_fatG}
\end{figure}

In Figure~\ref{fig:mean_curvature_Torus_fatG} we report on the time evolution of the minimum and maximum of mean curvature, i.e.~we plot $\min\{ H_h^n \}$ and $\max\{ H_h^n \}$. The curvature thresholds $H_i$, and the surgery time $T_1$ are marked again.

\begin{figure}[htbp]
	\centering
	\includegraphics[width=1\textwidth,clip,trim={110 155 50 180}]
	{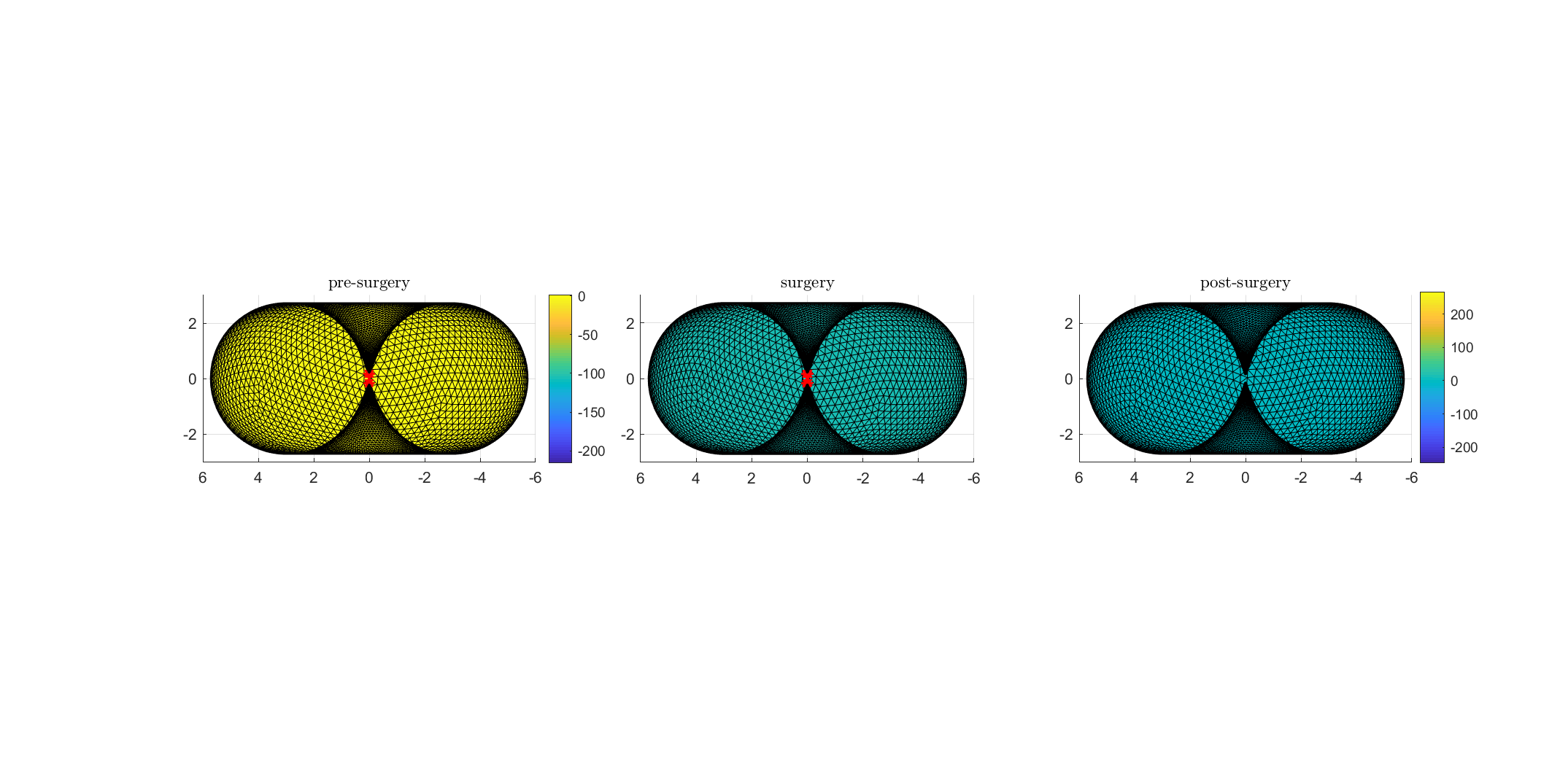}
	\includegraphics[width=1\textwidth,clip,trim={100 37 60 37}]
	{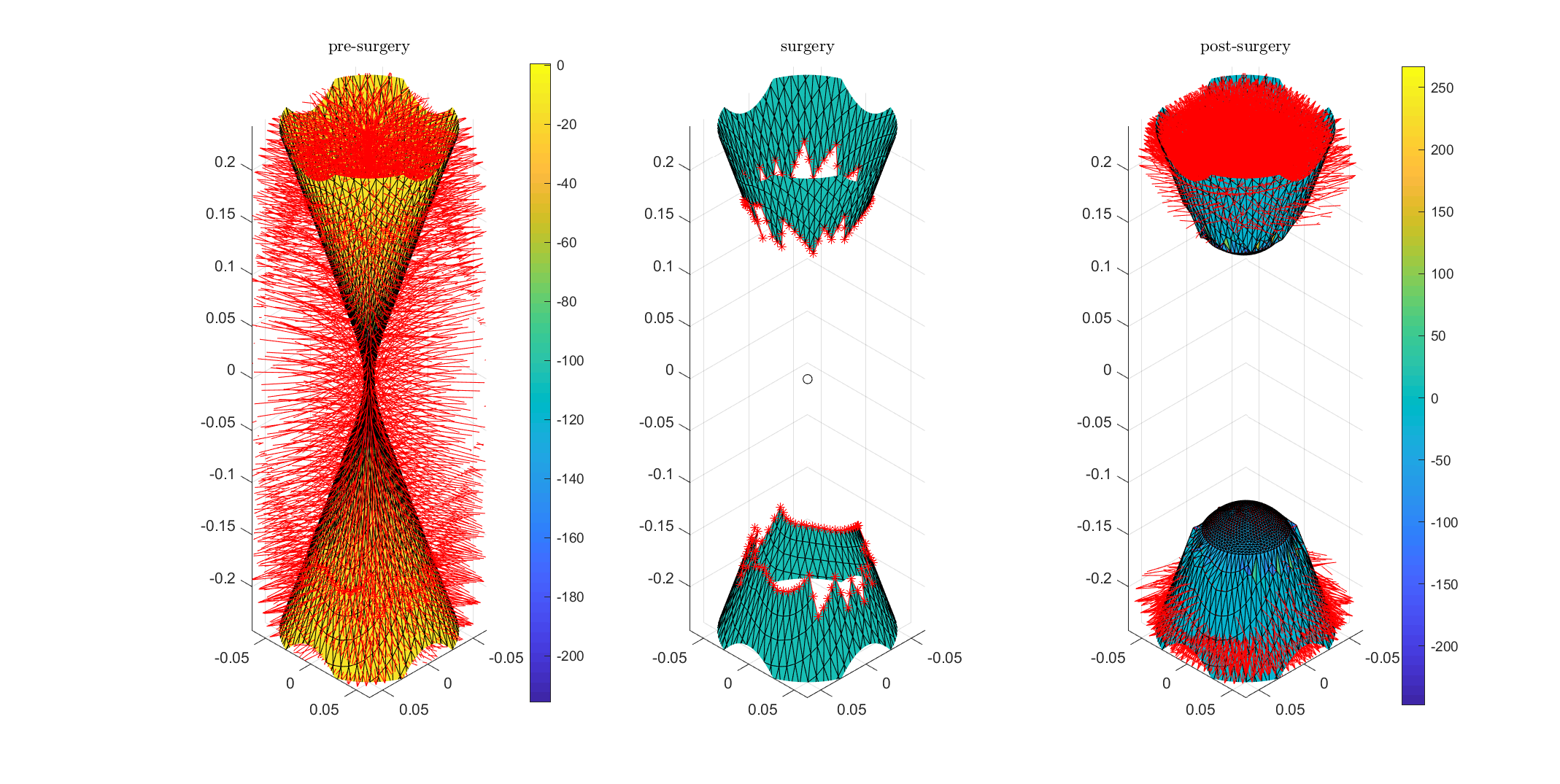}
	\caption{The steps of numerical surgery Algorithm~\ref{alg:numerical surgery} for the merging torus example, when $\max\{|H_h^n|\}$ exceeds $H_3$ at time $T_1 = 0.03433$. Note that in the close-up images the normal vectors point outwards, through the surface.}
	\label{fig:numerical_surgery_Torus_fatG}
\end{figure}

In Figure~\ref{fig:numerical_surgery_Torus_fatG} we again report on the numerical surgery for the merging torus example. At $T_1 = 0.03433$ an ``inverted'' pinch singularity is resolved using the extension described in Remark~\ref{remark:fat torus}.

\newpage

\section*{Acknowledgements}

The author thanks Christian Lubich and Ao Sun for our discussions on the topic. \bbk We also thank two referees who helped us to improve the mansucript. We are in particular grateful for the many references on singularities and their profiles. \ebk 

During the preparation of the manuscript and the codes the author was working at the Technical University of Munich and at the University of Regensburg.

The work of Bal\'azs Kov\'acs is funded by the Heisenberg Programme of the Deut\-sche For\-schungs\-gemeinschaft (DFG, German Research Foundation) -- \linebreak Project-ID 446431602, and by the DFG Graduiertenkolleg 2339 \emph{IntComSin} -- Project-ID 321821685.

\bibliographystyle{siamplain}
\bibliography{evolving_surface_literature}

\begin{thebibliography}{10}

\bibitem{ABG_MCFdiff_2}
{\sc H.~Abels, F.~B{\"u}rger, and H.~Garcke}, {\em Qualitative properties for a
  system coupling scaled mean curvature flow and diffusion}, arXiv:2205.02493,
  (2022).

\bibitem{ABG_MCFdiff_1}
{\sc H.~Abels, F.~B{\"u}rger, and H.~Garcke}, {\em Short time existence for
  coupling of scaled mean curvature flow and diffusion}, arXiv:2204.07626,
  (2022).

\bibitem{AlessandroniSinestrari_nhMCF}
{\sc R.~Alessandroni and C.~Sinestrari}, {\em Convexity estimates for a
  nonhomogeneous mean curvature flow}, Math. Z., 266 (2010), pp.~65--82,
  \url{https://doi.org/10.1007/s00209-009-0554-3}.

\bibitem{AltschulerAngenentGiga_1995}
{\sc S.~Altschuler, S.~B. Angenent, and Y.~Giga}, {\em Mean curvature flow
  through singularities for surfaces of rotation}, J. Geom. Anal., 5 (1995),
  pp.~293--358, \url{https://doi.org/10.1007/BF02921800},
  \url{https://doi.org/10.1007/BF02921800}.

\bibitem{AngenentVelazquez_1997}
{\sc S.~B. Angenent and J.~J.~L. Vel\'{a}zquez}, {\em Degenerate neckpinches in
  mean curvature flow}, J. Reine Angew. Math., 482 (1997), pp.~15--66,
  \url{https://doi.org/10.1515/crll.1997.482.15},
  \url{https://doi.org/10.1515/crll.1997.482.15}.

\bibitem{BaiLi_2022}
{\sc G.~Bai and B.~Li}, {\em Erratum: Convergence of {D}ziuk's semidiscrete
  finite element method for mean curvature flow of closed surfaces with
  high-order finite elements},  (2022).

\bibitem{BMPU_2012}
{\sc M.~Balazovjech, K.~Mikula, M.~Petr{\'a}{\v{s}}ov{\'a}, and J.~Urb{\'a}n},
  {\em Lagrangean method with topological changes for numerical modelling of
  forest fire propagation}, in Proceedings of ALGORITMY, 2012, pp.~42--52.

\bibitem{BaoZhao_2021}
{\sc W.~Bao and Q.~Zhao}, {\em A structure-preserving parametric finite element
  method for surface diffusion}, SIAM J. Numer. Anal., 59 (2021),
  pp.~2775--2799, \url{https://doi.org/10.1137/21M1406751}.

\bibitem{BarreiraElliottMadzvamuse2011}
{\sc R.~Barreira, C.~M. Elliott, and A.~Madzvamuse}, {\em The surface finite
  element method for pattern formation on evolving biological surfaces}, J.
  Math. Biol., 63 (2011), pp.~1095--1119,
  \url{https://doi.org/10.1007/s00285-011-0401-0}.

\bibitem{BGN2008}
{\sc J.~W. Barrett, H.~Garcke, and R.~N{\"u}rnberg}, {\em On the parametric
  finite element approximation of evolving hypersurfaces in {$\R^3$}}, J.
  Comput. Phys., 227 (2008), pp.~4281--4307.

\bibitem{BCH2006}
{\sc S.~Bartels, C.~Carstensen, and A.~Hecht}, {\em {$\rm P2Q2Iso2D=2D$}
  isoparametric {FEM} in {M}atlab}, J. Comput. Appl. Math., 192 (2006),
  pp.~219--250.

\bibitem{BRR_2018}
{\sc S.~Bartels, P.~Reiter, and J.~Riege}, {\em A simple scheme for the
  approximation of self-avoiding inextensible curves}, IMA J. Numer. Anal., 38
  (2018), pp.~543--565, \url{https://doi.org/10.1093/imanum/drx021}.

\bibitem{BG_1}
{\sc H.~Benninghoff and H.~Garcke}, {\em Efficient image segmentation and
  restoration using parametric curve evolution with junctions and topology
  changes}, SIAM J. Imaging Sci., 7 (2014), pp.~1451--1483,
  \url{https://doi.org/10.1137/130932430}.

\bibitem{BG_4}
{\sc H.~Benninghoff and H.~Garcke}, {\em Segmentation of three-dimensional
  images with parametric active surfaces and topology changes}, J. Sci.
  Comput., 72 (2017), pp.~1333--1367,
  \url{https://doi.org/10.1007/s10915-017-0401-3}.

\bibitem{MCF_generalised}
{\sc T.~Binz and B.~Kov{\'a}cs}, {\em A convergent finite element algorithm for
  generalized mean curvature flows of closed surfaces}, IMA J. Numer. Anal.,
  (2021).
\newblock doi:10.1093/imanum/drab043.

\bibitem{BrendleChoi_2019}
{\sc S.~Brendle and K.~Choi}, {\em Uniqueness of convex ancient solutions to
  mean curvature flow in {$\Bbb R^3$}}, Invent. Math., 217 (2019), pp.~35--76,
  \url{https://doi.org/10.1007/s00222-019-00859-4},
  \url{https://doi.org/10.1007/s00222-019-00859-4}.

\bibitem{BrendleChoi_2021}
{\sc S.~Brendle and K.~Choi}, {\em Uniqueness of convex ancient solutions to
  mean curvature flow in higher dimensions}, Geom. Topol., 25 (2021),
  pp.~2195--2234, \url{https://doi.org/10.2140/gt.2021.25.2195},
  \url{https://doi.org/10.2140/gt.2021.25.2195}.

\bibitem{BrendleHuisken2016}
{\sc S.~Brendle and G.~Huisken}, {\em Mean curvature flow with surgery of mean
  convex surfaces in {$\Bbb R^3$}}, Invent. Math., 203 (2016), pp.~615--654,
  \url{https://doi.org/10.1007/s00222-015-0599-3}.

\bibitem{BrendleHuisken2018}
{\sc S.~Brendle and G.~Huisken}, {\em Mean curvature flow with surgery of mean
  convex surfaces in three-manifolds}, J. Eur. Math. Soc. (JEMS), 20 (2018),
  pp.~2239--2257, \url{https://doi.org/10.4171/JEMS/811}.

\bibitem{diss_Buerger}
{\sc F.~B{\"u}rger}, {\em Interaction of mean curvature flow and diffusion},
  {PhD} thesis, University of Regensburg, Regensburg, Germany, 2021.
\newblock doi: 10.5283/epub.51215.

\bibitem{MCF_generic_I}
{\sc O.~Chodosh, K.~Choi, C.~Mantoulidis, and F.~Schulze}, {\em Mean curvature
  flow with generic initial data}, arXiv:2003.14344,  (2020).

\bibitem{MCF_generic_II}
{\sc O.~Chodosh, K.~Choi, and F.~Schulze}, {\em Mean curvature flow with
  generic initial data ii}, arXiv:2302.08409,  (2023).

\bibitem{CMP_survey}
{\sc T.~H. Colding, W.~P. Minicozzi, II, and E.~K. Pedersen}, {\em Mean
  curvature flow}, Bull. Amer. Math. Soc. (N.S.), 52 (2015), pp.~297--333,
  \url{https://doi.org/10.1090/S0273-0979-2015-01468-0}.

\bibitem{Deckelnick2000}
{\sc K.~Deckelnick}, {\em Error bounds for a difference scheme approximating
  viscosity solutions of mean curvature flow}, Interfaces Free Bound., 2
  (2000), pp.~117--142, \url{https://doi.org/10.4171/IFB/15}.

\bibitem{DeckelnickDziuk2003}
{\sc K.~Deckelnick and G.~Dziuk}, {\em A finite element level set method for
  anisotropic mean curvature flow with space dependent weight}, in Geometric
  analysis and nonlinear partial differential equations, Springer, Berlin,
  2003, pp.~249--264.

\bibitem{DeckelnickDE2005}
{\sc K.~Deckelnick, G.~Dziuk, and C.~M. Elliott}, {\em Computation of geometric
  partial differential equations and mean curvature flow}, Acta Numerica, 14
  (2005), pp.~139--232.

\bibitem{Demlow2009}
{\sc A.~Demlow}, {\em Higher--order finite element methods and pointwise error
  estimates for elliptic problems on surfaces}, SIAM J. Numer. Anal., 47
  (2009), pp.~805--807, \url{https://doi.org/10.1137/070708135}.

\bibitem{DemlowDziuk_2007}
{\sc A.~Demlow and G.~Dziuk}, {\em An adaptive finite element method for the
  {L}aplace-{B}eltrami operator on implicitly defined surfaces}, SIAM J. Numer.
  Anal., 45 (2007), pp.~421--442, \url{https://doi.org/10.1137/050642873}.

\bibitem{Dziuk90}
{\sc G.~Dziuk}, {\em An algorithm for evolutionary surfaces}, Numer. Math., 58
  (1990), pp.~603--611.

\bibitem{Ecker2012}
{\sc K.~Ecker}, {\em Regularity theory for mean curvature flow}, vol.~57 of
  Progress in Nonlinear Differential Equations and their Applications,
  Birkh\"{a}user Boston, Boston, MA, 2004.

\bibitem{EilksElliott2008}
{\sc C.~Eilks and C.~M. Elliott}, {\em Numerical simulation of dealloying by
  surface dissolution via the evolving surface finite element method}, J.
  Comput. Phys., 227 (2008), pp.~9727--9741,
  \url{https://doi.org/10.1016/j.jcp.2008.07.023}.

\bibitem{Elliott-Fritz-2017}
{\sc C.~M. Elliott and H.~Fritz}, {\em {On approximations of the curve
  shortening flow and of the mean curvature flow based on the DeTurck trick}},
  IMA J. Numer. Anal., 37 (2017), pp.~543--603.

\bibitem{MCFdiff}
{\sc C.~M. Elliott, H.~Garcke, and B.~Kov{\'a}cs}, {\em Numerical analysis for
  the interaction of mean curvature flow and diffusion on closed surfaces},
  Numer.~Math., 151 (2022), pp.~873--925,
  \url{https://doi.org/10.1007/s00211-022-01301-3}.

\bibitem{Espin_nhMCF}
{\sc T.~Espin}, {\em A pinching estimate for convex hypersurfaces evolving
  under a non-homogeneous variant of mean curvature flow}, Proc. Edinb. Math.
  Soc. (2), 65 (2022), pp.~376--391,
  \url{https://doi.org/10.1017/S001309152200013X}.

\bibitem{FengProhl2003}
{\sc X.~Feng and A.~Prohl}, {\em Numerical analysis of the {A}llen-{C}ahn
  equation and approximation for mean curvature flows}, Numer. Math., 94
  (2003), pp.~33--65, \url{https://doi.org/10.1007/s00211-002-0413-1}.

\bibitem{FierroPaolini1996}
{\sc F.~Fierro and M.~Paolini}, {\em Numerical evidence of fattening for the
  mean curvature flow}, Math. Models Methods Appl. Sci., 6 (1996),
  pp.~793--813, \url{https://doi.org/10.1142/S021820259600033X}.

\bibitem{FifeCahnElliott2001}
{\sc P.~C. Fife, J.~W. Cahn, and C.~M. Elliott}, {\em A free-boundary model for
  diffusion-induced grain boundary motion}, Interfaces Free Bound., 3 (2001),
  pp.~291--336, \url{https://doi.org/10.4171/IFB/42},
  \url{https://doi.org/10.4171/IFB/42}.

\bibitem{Gerhardt_pIMCF}
{\sc C.~Gerhardt}, {\em Non-scale-invariant inverse curvature flows in
  {E}uclidean space}, Calc. Var. Partial Differential Equations, 49 (2014),
  pp.~471--489, \url{https://doi.org/10.1007/s00526-012-0589-x}.

\bibitem{GonzalezMaddocks_1999}
{\sc O.~Gonzalez and J.~H. Maddocks}, {\em Global curvature, thickness, and the
  ideal shapes of knots}, Proc. Natl. Acad. Sci. USA, 96 (1999),
  pp.~4769--4773, \url{https://doi.org/10.1073/pnas.96.9.4769}.

\bibitem{Hamilton95}
{\sc R.~S. Hamilton}, {\em The formation of singularities in the {R}icci flow},
  in Surveys in differential geometry, {V}ol. {II} ({C}ambridge, {MA}, 1993),
  Int. Press, Cambridge, MA, 1995, pp.~7--136.

\bibitem{Hamilton97}
{\sc R.~S. Hamilton}, {\em Four-manifolds with positive isotropic curvature},
  Comm. Anal. Geom., 5 (1997), pp.~1--92,
  \url{https://doi.org/10.4310/CAG.1997.v5.n1.a1}.

\bibitem{Haslhofer_2015}
{\sc R.~Haslhofer}, {\em Uniqueness of the bowl soliton}, Geom. Topol., 19
  (2015), pp.~2393--2406.

\bibitem{HaslhoferKleiner_2017_survey}
{\sc R.~Haslhofer and B.~Kleiner}, {\em Mean curvature flow of mean convex
  hypersurfaces}, Comm. Pure Appl. Math., 70 (2017), pp.~511--546,
  \url{https://doi.org/10.1002/cpa.21650},
  \url{https://doi.org/10.1002/cpa.21650}.

\bibitem{HK_2017}
{\sc R.~Haslhofer and B.~Kleiner}, {\em Mean curvature flow with surgery}, Duke
  Math. J., 166 (2017), pp.~1591--1626,
  \url{https://doi.org/10.1215/00127094-0000008X}.

\bibitem{HuLi_2022}
{\sc J.~Hu and B.~Li}, {\em Evolving finite element methods with an artificial
  tangential velocity for mean curvature flow and {W}illmore flow}, Numer.
  Math., 152 (2022), pp.~127--181.

\bibitem{Huisken1984}
{\sc G.~Huisken}, {\em Flow by mean curvature of convex surfaces into spheres},
  J. Differential Geometry, 20 (1984), pp.~237--266.

\bibitem{Huisken1990}
{\sc G.~Huisken}, {\em Asymptotic behavior for singularities of the mean
  curvature flow}, J. Differential Geom., 31 (1990), pp.~285--299,
  \url{http://projecteuclid.org/euclid.jdg/1214444099}.

\bibitem{HuiskenIlmanen}
{\sc G.~Huisken and T.~Ilmanen}, {\em The inverse mean curvature flow and the
  {R}iemannian {P}enrose inequality}, J. Differential Geom., 59 (2001),
  pp.~353--437.

\bibitem{HuiskenPolden}
{\sc G.~Huisken and A.~Polden}, {\em Geometric evolution equations for
  hypersurfaces}, in Calculus of variations and geometric evolution problems
  ({C}etraro, 1996), vol.~1713 of Lecture Notes in Math., Springer, Berlin,
  1999, pp.~45--84.

\bibitem{HuiskenSinestrari_1999a}
{\sc G.~Huisken and C.~Sinestrari}, {\em Convexity estimates for mean curvature
  flow and singularities of mean convex surfaces}, Acta Math., 183 (1999),
  pp.~45--70, \url{https://doi.org/10.1007/BF02392946}.

\bibitem{HuiskenSinestrari_1999b}
{\sc G.~Huisken and C.~Sinestrari}, {\em Mean curvature flow singularities for
  mean convex surfaces}, Calc. Var. Partial Differential Equations, 8 (1999),
  pp.~1--14, \url{https://doi.org/10.1007/s005260050113}.

\bibitem{HuiskenSinestrari}
{\sc G.~Huisken and C.~Sinestrari}, {\em Mean curvature flow with surgeries of
  two-convex hypersurfaces}, Invent. Math., 175 (2009), pp.~137--221,
  \url{https://doi.org/10.1007/s00222-008-0148-4}.

\bibitem{highorderESFEM}
{\sc B.~Kov\'{a}cs}, {\em High-order evolving surface finite element method for
  parabolic problems on evolving surfaces}, IMA J. Numer. Anal., 38 (2018),
  pp.~430--459.

\bibitem{MCF}
{\sc B.~Kov{\'a}cs, B.~Li, and C.~Lubich}, {\em A convergent evolving finite
  element algorithm for mean curvature flow of closed surfaces}, Numer. Math.,
  143 (2019), pp.~797--853.

\bibitem{MCF_soldriven}
{\sc B.~Kov{\'a}cs, B.~Li, and C.~Lubich}, {\em A convergent algorithm for
  forced mean curvature flow driven by diffusion on the surfaces}, Interfaces
  Free Bound., 22 (2020), pp.~443--464.

\bibitem{Willmore}
{\sc B.~Kov{\'a}cs, B.~Li, and C.~Lubich}, {\em A convergent evolving finite
  element algorithm for {W}illmore flow of closed surfaces}, Numer.~Math., 149
  (2021), pp.~595--643.

\bibitem{Lantelme}
{\sc M.~Lantelme}, {\em A posteriori error estimates and adaptive algorithms
  for parabolic surface {PDEs}}, master thesis, Technical University of Munich,
  Munich, German, 2023.

\bibitem{Li_2021}
{\sc B.~Li}, {\em Convergence of {D}ziuk's semidiscrete finite element method
  for mean curvature flow of closed surfaces with high-order finite elements},
  SIAM J. Numer. Anal., 59 (2021), pp.~1592--1617,
  \url{https://doi.org/10.1137/20M136935X}.

\bibitem{Mantegazza}
{\sc C.~Mantegazza}, {\em {Lecture Notes on Mean Curvature Flow}}, Progress in
  Mathematics, Volume 290. Birkh\"auser, Corrected Printing 2012.

\bibitem{MBO}
{\sc B.~Merriman, J.~K. Bence, and S.~Osher}, {\em Diffusion generated motion
  by mean curvature}, Department of Mathematics, University of California, Los
  Angeles, 1992.

\bibitem{MikulaUrban_2012}
{\sc K.~Mikula and J.~Urb{\'a}n}, {\em New fast and stable lagrangean method
  for image segmentation}, in 2012 5th International Congress on Image and
  Signal Processing, IEEE, 2012, pp.~688--696.

\bibitem{Mitchell_2}
{\sc I.~M. Mitchell}, {\em A toolbox of level set methods (version 1.1)},
  Department of Computer Science, University of British Columbia, Vancover, BC,
  Canada,  (2007),
  \url{https://www.cs.ubc.ca/~mitchell/ToolboxLS/toolboxLS-1.1.pdf}.
\newblock Tech.~Rep.~TR-2007-11 (online).

\bibitem{Mitchell_1}
{\sc I.~M. Mitchell}, {\em The flexible, extensible and efficient toolbox of
  level set methods}, J. Sci. Comput., 35 (2008), pp.~300--329,
  \url{https://doi.org/10.1007/s10915-007-9174-4}.

\bibitem{Nochetto_etal_1994}
{\sc R.~H. Nochetto, M.~Paolini, and C.~Verdi}, {\em Optimal interface error
  estimates for the mean curvature flow}, Ann. Scuola Norm. Sup. Pisa Cl. Sci.
  (4), 21 (1994), pp.~193--212.

\bibitem{Nochetto_etal_1996}
{\sc R.~H. Nochetto and C.~Verdi}, {\em Combined effect of explicit
  time-stepping and quadrature for curvature driven flows}, Numer. Math., 74
  (1996), pp.~105--136, \url{https://doi.org/10.1007/s002110050210}.

\bibitem{OOTT_2011}
{\sc A.~Oberman, S.~Osher, R.~Takei, and Y.-H.~R. Tsai}, {\em Numerical methods
  for anisotropic mean curvature flow based on a discrete time variational
  formulation}, Commun. Math. Sci., 9 (2011), pp.~637--662,
  \url{https://doi.org/10.4310/CMS.2011.v9.n3.a1}.

\bibitem{Olshanskii_etal_2021}
{\sc M.~Olshanskii, X.~Xu, and V.~Yushutin}, {\em A finite element method for
  {A}llen--{C}ahn equation on deforming surface}, Comput. Math. Appl., 90
  (2021), pp.~148--158, \url{https://doi.org/10.1016/j.camwa.2021.03.018}.

\bibitem{OsherSethian_1988}
{\sc S.~Osher and J.~A. Sethian}, {\em Fronts propagating with
  curvature-dependent speed: algorithms based on {H}amilton-{J}acobi
  formulations}, J. Comput. Phys., 79 (1988), pp.~12--49,
  \url{https://doi.org/10.1016/0021-9991(88)90002-2}.

\bibitem{Perelman_1}
{\sc G.~Perelman}, {\em The entropy formula for the {R}icci flow and its
  geometric applications},  (2002).
\newblock arXiv:math/0211159.

\bibitem{Perelman_3}
{\sc G.~Perelman}, {\em Finite extinction time for solutions to the {R}icci
  flow on certain three-manifolds},  (2003).
\newblock arXiv:math/0307245.

\bibitem{Perelman_2}
{\sc G.~Perelman}, {\em Ricci flow with surgery on three-manifolds},  (2003).
\newblock arXiv:math/0303109.

\bibitem{distmesh}
{\sc P.-O. Persson and G.~Strang}, {\em A simple mesh generator in {MATLAB}},
  SIAM Review, 46 (2004), pp.~329--345.

\bibitem{Scheuer_pIMCF}
{\sc J.~Scheuer}, {\em Pinching and asymptotical roundness for inverse
  curvature flows in {E}uclidean space}, J. Geom. Anal., 26 (2016),
  pp.~2265--2281, \url{https://doi.org/10.1007/s12220-015-9627-1}.

\bibitem{Schulze_diss}
{\sc F.~Schulze}, {\em Nichtlineare Evolution von Hyperfl{\"a}chen entlang
  ihrer mittleren Kr{\"u}mmung}, {PhD} thesis, University of T\"ubingen,
  T\"ubingen, Germany, 2002.

\bibitem{Schulze_1}
{\sc F.~Schulze}, {\em Evolution of convex hypersurfaces by powers of the mean
  curvature}, Math. Z., 251 (2005), pp.~721--733,
  \url{https://doi.org/10.1007/s00209-004-0721-5}.

\bibitem{Sinestrari_survey}
{\sc C.~Sinestrari}, {\em Singularities of mean curvature flow and flow with
  surgeries}, in Surveys in differential geometry. {V}ol. {XII}. {G}eometric
  flows, vol.~12 of Surv. Differ. Geom., Int. Press, Somerville, MA, 2008,
  pp.~303--332, \url{https://doi.org/10.4310/SDG.2007.v12.n1.a8}.

\bibitem{White_2000}
{\sc B.~White}, {\em The size of the singular set in mean curvature flow of
  mean-convex sets}, J. Amer. Math. Soc., 13 (2000), pp.~665--695,
  \url{https://doi.org/10.1090/S0894-0347-00-00338-6}.

\bibitem{White_2003}
{\sc B.~White}, {\em The nature of singularities in mean curvature flow of
  mean-convex sets}, J. Amer. Math. Soc., 16 (2003), pp.~123--138,
  \url{https://doi.org/10.1090/S0894-0347-02-00406-X}.

\bibitem{WittwerAland2022}
{\sc L.~D. Wittwer and S.~Aland}, {\em A computational model of self-organized
  shape dynamics of active surfaces in fluids}, arXiv:2203.00099,  (2022).

\bibitem{BaoJiangZhao_2020}
{\sc Q.~Zhao, W.~Jiang, and W.~Bao}, {\em A parametric finite element method
  for solid-state dewetting problems in three dimensions}, SIAM J. Sci.
  Comput., 42 (2020), pp.~B327--B352, \url{https://doi.org/10.1137/19M1281666}.

\end{thebibliography}
\end{document}